\documentclass[11pt,reqno]{amsart}

\usepackage[margin=1.08in]{geometry}
\usepackage{amsmath,amssymb,mathtools,bm,mathrsfs}
\usepackage{enumitem}
\usepackage[protrusion=true,expansion=false]{microtype}
\usepackage{xcolor}
\usepackage[colorlinks=true,linkcolor=blue!55!black,citecolor=blue!55!black,
            urlcolor=blue!60!black]{hyperref}

\allowdisplaybreaks
\numberwithin{equation}{section}

\newtheorem{theorem}{Theorem}[section]
\newtheorem{proposition}[theorem]{Proposition}
\newtheorem{lemma}[theorem]{Lemma}

\newtheorem{assumption}[theorem]{Assumption}
\theoremstyle{definition}

\newtheorem{example}[theorem]{Example}
\theoremstyle{remark}
\newtheorem{remark}[theorem]{Remark}

\newcommand{\R}{\mathbb R}
\newcommand{\N}{\mathbb N}
\newcommand{\dd}{\,\mathrm d}

\newcommand{\av}[1]{\left\langle #1\right\rangle}
\newcommand{\calF}{\mathcal F}
\newcommand{\calG}{\mathcal G}
\newcommand{\calL}{\mathcal L}
\newcommand{\calS}{\mathcal S}

\newcommand{\calX}{\mathscr X}
\newcommand{\calY}{\mathscr Y}
\newcommand{\Ran}{\operatorname{Ran}}
\newcommand{\Ker}{\operatorname{Ker}}

\title[Local and global bifurcation for stationary RVM]
{Local and Global Bifurcation\\
for Axisymmetric Stationary Relativistic Vlasov--Maxwell Equilibria}

\author{Katherine Zhiyuan Zhang}
\address{Department of Mathematics, Northeastern University}
\dedicatory{Dedicated to Walter A. Strauss}

\begin{document}
\raggedbottom

\begin{abstract}
We study bifurcation of stationary two-species relativistic
Vlasov--Maxwell equilibria on bounded axisymmetric domains.
Velocity-moment and elliptic spectral bounds give
computable tests for simple crossings, their number and separation,
while the local theorem gives a branch-direction formula.  Explicit
two-species families exhibit
electric-only, magnetic-only, and joint field-norm escape. 
\end{abstract}

\maketitle
\enlargethispage{8pt}
\tableofcontents

\section{Introduction}
\label{sec:introduction}

Stationary relativistic Vlasov--Maxwell (RVM) equilibria couple the
collisionless particle distribution to the field it generates.  Their
existence, secondary bifurcation, and dynamical stability are distinct
questions.  A regular continuation family need not bifurcate; conversely,
a zero mode of the stationary field equations need not generate a
growing mode of the time-dependent kinetic system.  The aim here is to
make the stationary bifurcation hypotheses computable from distribution
profiles and domain spectra, and to prove norm escape in selected global
families.  The assertions below concern stationary equilibria and their
stationary linearizations; they do not supply a growing-mode or
spectral-stability theorem for the time-dependent RVM system.

We start from the bounded three-dimensional axisymmetric setting of
Zhang \cite{Zhang2023}, with specular particle reflection and perfectly
conducting field boundary conditions.  For
\[
 E=-\nabla\phi,\qquad
 \bm A=A(r,z)e_\varphi,\qquad B=\nabla\times\bm A,
\]
the two-species invariant ansatz is
\begin{equation}
 f_K^\sigma(x,v)
 =\mu_K^\sigma(e_\sigma(x,v),p_\sigma(x,v)),
 \qquad \sigma\in\{+1,-1\},                      \label{eq:intro-ansatz}
\end{equation}
where \(e_\sigma\) and \(p_\sigma\) denote relativistic particle
energy and canonical angular momentum.  The kinetic equations are then
satisfied along characteristics, while Maxwell's equations reduce to a
coupled Dirichlet problem for \((A,\phi)\).  Theorem~3.2 in \cite{Zhang2023} already
provides a two-species global continuum with a loop alternative, and
Theorem~3.3 obtains an unbounded locally analytic curve under a strict
single-species assumption \cite{Zhang2023}.  The same paper also discusses
parameter-dependent stability changes; we do not claim that phenomenon
itself is new.  Those curves start at a
regular base point under the paper's hypotheses; the present questions
instead concern a singular point of the stationary derivative along a
reference family and the secondary solutions emerging there.

\subsection*{Previous results and comparison}

Stationary collisionless plasmas in bounded domains and relativistic
special-form equilibria have a substantial history
\cite{Rein1992,Poupaud,Degond1990,BattFabian,Braasch,
Bostan2009,Schaeffer2005,Suzuki2008}; Weber's confined relativistic
states are constructed in an infinite cylinder \cite{Weber2020}.
Classical stationary Vlasov--Maxwell bifurcation was studied for
two-particle and multicomponent distributions by Rudykh, Sidorov,
Sinitsyn, and collaborators
\cite{RudykhSidorovSinitsyn1989,SidorovSinitsyn1996,
SidorovSinitsyn1997,SidorovSinitsyn2000,RendonEtAl2016}.
In particular, \cite{SidorovSinitsyn1996} proves an unbounded
first-mode component in its classical problem, and
\cite{SidorovSinitsyn1997} analyzes bifurcation points and nontrivial
branch asymptotics.  Steady classical equilibria and stability are
also treated in \cite{MarkovEtAl1992}.  More recent examples include
pinches with additional cylindrical/vertical symmetry whose principal
construction is classical but whose Section~9 also treats relativistic
adaptations \cite{HolmerZhang2026}.  These results preclude any claim
that multispecies branches or unbounded stationary VM curves first
arise here.  They use different kinetic energies, geometric reductions,
boundary conditions, or profile classes from the bounded axisymmetric
RVM problem considered below.
Guo and
Strauss gave a homogeneous equilibrium stable to purely electric
perturbations but unstable to certain electromagnetic ones
\cite{GuoStrauss2000}.  For relativistic plasmas, Lin and Strauss
developed stability operators that retain the effect of transport along
particle trajectories \cite{LinStrauss2007,LinStrauss2008}.
Nguyen and Strauss treated specular reflection first in a problem
with radial symmetry \cite{NguyenStrauss2013} and then in the solid
torus \cite{NguyenStrauss2014}; Zhang treated general bounded
axisymmetric domains \cite{Zhang2019}.

Canonical-momentum-dependent profiles have antecedents in current
sheets and rigid-rotor models
\cite{Harris1962,Channell1976,ChenDavidsonPakter1997}; the thermal
factor is related to the Maxwell--J\"uttner distribution
\cite{Juttner1911}.  Our explicit laws are mathematically chosen
representatives of these modelling classes, not measured fits.
Numerical plasma studies also use ``bifurcation'' for transitions
between computed profile regimes, including cylindrical equilibria
and diode configurations \cite{CordonnierEtAl2022,
SidorovEtAlDiode2026}.  Here it means a singular stationary
Fr\'echet derivative and, where asserted, a simple kernel with a
verified transverse crossing.

There is also a connection to Strauss's work on local/global continuation in steady water waves \cite{ConstantinStraussVarvaruca2016} and rotating gravitational Vlasov-Poisson equilibria \cite{StraussWu2017}.
For rotating Vlasov--Poisson galaxies, Strauss and Wu subsequently
proved a global alternative in physical terms: the spatial supports
are unbounded, or both rotation speeds and maximum densities are
unbounded \cite{StraussWu2023}.

Our abstract input is the classical local simple-eigenvalue theorem
\cite{CrandallRabinowitz1971} and global degree, positive-cone,
and analytic-continuation methods
\cite{Rabinowitz1971,DancerPositive1973,DancerAnalytic1973,
DancerGlobal1973,BuffoniToland}.  Neither block elimination nor a
self-adjoint Schur complement is new as an operator device: related
operators occur in axisymmetric RVM stability analysis
\cite{NguyenStrauss2014,Zhang2019}.  Those dynamical operators include
transport-orbit projections and answer a different question from
singularity of the stationary Maxwell map.  Likewise, our resonant
norm-escape and sublinear scalar arguments specialize established
theory \cite{RabinowitzInfinity1973,BrezisOswald1986,Klimsiak2026},
and do not by themselves give a general classification of global components.

\subsection*{Contributions and technical approach}

Relative to the cited literature, the contribution is the combination
of explicit stationary RVM tests and realizations in this geometry.
First, an integration-by-parts identity in the two-species moment
derivatives produces a self-adjoint magnetic Schur complement; its
simple zero and crossing form distinguish genuine secondary
bifurcation from regular continuation.  Electric--magnetic coupling and
opposite species charges make scalar sign arguments unreliable.  We
therefore derive electric coercivity, endpoint Rayleigh inequalities,
quantitative monotonicity, and a second-mode gap directly from velocity
moments.  Fixed magnetic-operator eigenvalue enclosures then locate,
count, and separate crossings without assuming knowledge of the
parameter-dependent Schur eigenvalues.  General even--odd two-species
profiles on variable solid tori realize both diagonal and genuinely
coupled crossings.  The neutral local pitchfork mechanism is also
realized globally by the monotone-curvature class in Section
\ref{sec:global-bif}.

Second, the global alternative alone permits norm escape, return to
another bifurcation point, or approach to the admissibility boundary;
it does not identify an escaping field.  Section~\ref{sec:global-bif}
instead constructs explicit two-species equilibrium curves and uses
principal-eigenfunction tests, finite-amplitude bounds, scalar
comparison, and positive-cone continuation to distinguish electric-only,
magnetic-only, and joint norm escape.  The electric-only graph is a
regular continuation, not a secondary bifurcation curve.  The magnetic
examples include a finite-parameter asymptotic resonance and a
counter-rotating Maxwell--J\"uttner pitchfork whose positive component
is confined to \(0<K<K_*\) but has unbounded magnetic norm.  These
examples do not classify all components allowed by the abstract
alternative or imply dynamical instability.  The underlying degree,
sublinear uniqueness, and positive-cone principles are classical;
the contribution is their explicit verification for the stated
two-species relativistic profiles and boundary conditions.


\subsection*{Organization}

Section \ref{sec:setup} gives the stationary reduction and function
spaces.  Section \ref{sec:linearization} proves the Schur reduction,
local branch statements, and the global alternative.
Section \ref{sec:local-bif} develops computable local crossing,
counting, and gap tests with variable-torus and
general-profile examples.  Section \ref{sec:global-bif} gives explicit
two-species norm-escape curves and finite tests for electric, magnetic,
and joint field growth.  

\subsection*{Disclosure of AI assistance}
During the preparation of this manuscript, the author provided the mathematical intuition, strategic directions, and conceptual interpretations, that is, the bifurcation and norm growth criteria in the manuscript. OpenAI's GPT-6 is used to assist with numerical computation and carrying out auxiliary calculations when identifying the concrete examples. 




\section{Stationary reduction and functional setting}
\label{sec:setup}

\subsection{The RVM system and its invariants}

Let $\Omega\subset\R^3$ be a connected, bounded, axisymmetric domain and
write cylindrical coordinates as $(r,\varphi,z)$.  For clarity of elliptic
regularity and the Hopf lemma we assume throughout that
\begin{equation}
 \partial\Omega\in C^{2,\alpha},\qquad 0<\alpha<1.
 \label{eq:domain}
\end{equation}

We consider the 3D relativistic Vlasov-Maxwell (RVM) system on an axisymmetric bounded domain $\Omega \subset \mathbb{R}^3$. We study a plasma with two species (ion and electron) with non-negative distribution functions $f^\pm(t, x, v)$, where $t \geq 0$, $x \in \Omega \subset \mathbb{R}^3$, $v \in \mathbb{R}^3$.
The system RVM is
\begin{equation} \label{Vlasovtwospecies}
\partial_t f^{\pm} + \hat{v} \cdot \nabla_x f^{\pm} \pm (E + \hat{v} \times B   ) \cdot \nabla_v f^{\pm} =0 ,
\end{equation}
\begin{equation}  \label{Maxwelltwospecies-1}
\nabla_x \cdot E = \rho = \int_{\mathbb{R}^3} (f^+ - f^-)dv, \  \nabla_x \cdot B =0 ,
\end{equation}
\begin{equation}  \label{Maxwelltwospecies-2}
\partial_t E -\nabla_x \times B = -\textbf{j} = -  \int_{\mathbb{R}^3} \hat{v} (f^+ - f^-)dv, \  \partial_t  B + \nabla_x \times E =0 .
\end{equation}
In this system, $f^\pm (t, x, v) \geq 0$ is the density distribution of ions ($+$) and electrons ($-$). $x \in \Omega$ is the particle position. $v \in \mathbb{R}^3$ is the particle momentum, $\langle v\rangle = \sqrt{1+ v^2}$ is the particle kinetic energy, and $\hat{v} = v/ \langle v \rangle$ is the particle velocity. Also, $E$ is the electric field, $B$ is the internal magnetic field, so $ \pm (E + \hat{v} \times B ) $ is the total (internal) electromagnetic force. Moreover, the charge density $\rho$ and the current density $\textbf{j}$ are defined as
\begin{equation}
\rho = \int_{\mathbb{R}^3} (f^+ - f^-) dv , \qquad \textbf{j} =  \int_{\mathbb{R}^3} \hat{v} (f^+ - f^-) dv  .
\end{equation}

For species sign $\sigma\in\{+1,-1\}$, the stationary characteristics are
generated by
\[
 \dot X=\widehat V,\qquad
 \dot V=\sigma\bigl(E(X)+\widehat V\times B(X)\bigr),
 \qquad \widehat V=\frac{V}{\sqrt{1+|V|^2}}.
\]
We consider the case
\[
 E=-\nabla\phi,\qquad
 \bm A=A(r,z)e_\varphi,\qquad B=\nabla\times\bm A,
\]
the particle energy and canonical angular momentum are
\begin{equation}
 e_\sigma=\av v+\sigma\phi,\qquad
 p_\sigma=r(v_\varphi+\sigma A),\qquad
 \av v=\sqrt{1+|v|^2}.
 \label{eq:invariants}
\end{equation}
They are conserved between reflections.  Specular reflection changes only
the normal component of $v$ and hence preserves both quantities.  It
follows that
\begin{equation}
 f_K^\sigma(x,v)=\mu_K^\sigma(e_\sigma(x,v),p_\sigma(x,v))
 \label{eq:ansatz}
\end{equation}
solves the stationary Vlasov equation and the specular boundary condition
whenever the moments below are finite.
Throughout the physical parameter range $K\ge0$ we require
$\mu_K^\sigma\ge0$ and finiteness of these moments.  A signed extension
to $K<0$, when used below, is only an auxiliary analytic device and is
not interpreted as a particle density.

Define the charge and azimuthal-current moment functions
\begin{align}
 \rho_K(r,a,b)
 &=\int_{\R^3}\!\left[
 \mu_K^+(\av v+a,r(v_\varphi+b))
 -\mu_K^-(\av v-a,r(v_\varphi-b))\right]\dd v,
 \label{eq:rho-moment}\\
 j_K(r,a,b)
 &=\int_{\R^3}\!\widehat v_\varphi\left[
 \mu_K^+(\av v+a,r(v_\varphi+b))
 -\mu_K^-(\av v-a,r(v_\varphi-b))\right]\dd v.
 \label{eq:j-moment}
\end{align}
We work in the fixed homogeneous potential, grounding, and magnetic-flux
class
\begin{equation}
 A|_{\partial\Omega}=\phi|_{\partial\Omega}=0.
 \label{eq:fixed-potential-class}
\end{equation}
This is a choice within, not a consequence of, the perfect-conductor
conditions $E_{\rm tan}=0$ and $B_n=0$.  On a multiply connected
domain the latter permit boundary constants and harmonic fields; for
example $A=c/r$ on an annular cylinder has vanishing curl.  
Maxwell's equations reduce to
\begin{equation}
 \left\{
 \begin{aligned}
 \left(-\Delta+\frac1{r^2}\right)A&=j_K(r,\phi,A),\\
 -\Delta\phi&=\rho_K(r,\phi,A),\\
 A|_{\partial\Omega}&=\phi|_{\partial\Omega}=0.
 \end{aligned}\right.
 \label{eq:elliptic-system}
\end{equation}
Once \eqref{eq:elliptic-system} is solved, \eqref{eq:ansatz} and the fields
$E=-\nabla\phi$, $B=\nabla\times(Ae_\varphi)$ give a stationary RVM
equilibrium.

\subsection{Spaces and nonlinear operators}

Fix $p>3$ and let the subscript ``ax'' denote independence of $\varphi$.
Set
\begin{align}
 X_\phi&=W^{2,p}_{\rm ax}(\Omega)\cap W^{1,p}_0(\Omega),\label{eq:Xphi}\\
 X_A&=\left\{A\in W^{1,p}_{0,\rm ax}(\Omega):
              e^{i\varphi}A\in W^{2,p}(\Omega;\mathbb C)\right\},
 \label{eq:XA}\\
 \calX&=X_A\times X_\phi,\qquad
 \calY=L^p_{\rm ax}(\Omega)\times L^p_{\rm ax}(\Omega).
\end{align}
The norm on $X_A$ is $\|e^{i\varphi}A\|_{W^{2,p}}$.  This encodes the
apparent singularity at $r=0$, since
\begin{equation}
 -\Delta(e^{i\varphi}A)
 =e^{i\varphi}\left(-\Delta+\frac1{r^2}\right)A.
 \label{eq:toroidal-identity}
\end{equation}
Write
\[
 L_A=-\Delta+r^{-2},\qquad L_\phi=-\Delta
\]
with homogeneous Dirichlet data.  Standard elliptic theory and
\eqref{eq:toroidal-identity} give bounded inverses
\[
 L_A^{-1}:L^p_{\rm ax}\to X_A,\qquad
 L_\phi^{-1}:L^p_{\rm ax}\to X_\phi.
\]

We use the strong residual
\begin{equation}
 \calF(A,\phi,K)=
 \binom{L_AA-j_K(r,\phi,A)}{L_\phi\phi-\rho_K(r,\phi,A)}
 \in\calY
 \label{eq:strong-F}
\end{equation}
and the equivalent fixed-point residual
\begin{equation}
 \calG(A,\phi,K)=
 \binom{A-L_A^{-1}j_K(r,\phi,A)}
       {\phi-L_\phi^{-1}\rho_K(r,\phi,A)}
 \in\calX.
 \label{eq:fixed-G}
\end{equation}
This formulation removes an ambiguity in the spaces used in the proof of
Theorem 3.2 of \cite{Zhang2023}: the strong map naturally goes from
$W^{2,p}$ to $L^p$, while the inverse-elliptic map acts on one Banach
space and is identity minus compact.

We impose the following standard regularity hypothesis.  It is a concise
version of the decay and factorial derivative estimates in Assumptions
3.1 and (3.21) of \cite{Zhang2023}.

\begin{assumption}[Analytic moment hypothesis]
\label{ass:analytic}
On every bounded set of $(a,b,K)$, the profiles and all derivatives used
below have a velocity-integrable majorant.  The functions
$(r,a,b,K)\mapsto(\rho_K,j_K)$ are real analytic in $(a,b,K)$ and
continuous in $r\in[0,d]$, where $d=\sup_{x\in\Omega}r(x)$.  More
precisely, the associated moment Nemytskii maps are real analytic from
$\calX\times\R$ to $L^p\times L^p$ and take bounded sets to bounded sets.
\end{assumption}

For example, locally uniform factorial bounds of the form
\[
 |\partial_K^\ell\partial_e^m\partial_p^n\mu_K^\sigma(e,p)|
 \le C^{\ell+m+n+1}\ell!m!n!(1+|e|)^{-\delta},\qquad \delta>3,
\]
on bounded $K$-intervals, with a velocity-integrable majorant after the
substitution \eqref{eq:ansatz}, imply this Banach-space analyticity.  Under
Assumption \ref{ass:analytic}, $\calF$ and $\calG$ are real analytic, and
\begin{equation}
 (z,K)\longmapsto\calG(z,K)-z
 \quad\hbox{is compact on bounded subsets of $\calX\times\R$.}
 \label{eq:compact}
\end{equation}
Indeed, first extract a convergent parameter subsequence from any bounded
sequence $(z_n,K_n)$.  The fields
$e^{i\varphi}A_n$ and the potentials $\phi_n$ have subsequences converging
in $C^1$ because $W^{2,p}\Subset C^1$ for $p>3$.  In particular
$A_n\to A$ and $\phi_n\to\phi$ uniformly (no scalar $C^1$ claim for $A$
at the axis is needed).  Continuity and local domination of the moment
maps give convergence of the sources in $L^p$, and bounded elliptic
inversion gives convergence in $\calX$.  In particular,
$D_{(A,\phi)}\calG$ is Fredholm of index zero at every point.

\section{Linearization and abstract bifurcation theory}
\label{sec:linearization}

\subsection{Linearization and Schur reduction}

Let $z_0(K)=(A_0(K),\phi_0(K))$ be a reference equilibrium and put
\begin{equation}
 E_\sigma=\av v+\sigma\phi_0,\qquad
 P_\sigma=r(v_\varphi+\sigma A_0).
 \label{eq:EP-ref}
\end{equation}
All profile derivatives below are evaluated at the reference states in
\eqref{eq:EP-ref}.  Define the coefficient functions
\begin{align}
 c_{11}&=r\int_{\R^3}\widehat v_\varphi
       (\mu^+_{K,p}+\mu^-_{K,p})\dd v,
 &c_{12}&=\int_{\R^3}\widehat v_\varphi
       (\mu^+_{K,e}+\mu^-_{K,e})\dd v,\label{eq:c11c12}\\
 c_{21}&=r\int_{\R^3}
       (\mu^+_{K,p}+\mu^-_{K,p})\dd v,
 &c_{22}&=\int_{\R^3}
       (\mu^+_{K,e}+\mu^-_{K,e})\dd v.\label{eq:c21c22}
\end{align}
For $h=(a,b)$, direct differentiation of \eqref{eq:strong-F} gives
\begin{equation}
 \calL(K)h:=D_z\calF(z_0(K),K)h
 =\begin{pmatrix}
 L_A-c_{11}&-c_{12}\\
 -c_{21}&L_\phi-c_{22}
 \end{pmatrix}\binom ab.
 \label{eq:full-linearization}
\end{equation}

\begin{lemma}[Off-diagonal cancellation]
\label{lem:cancellation}
Suppose, for each fixed $(x,K)$ under consideration,
$\mu_K^\sigma(E_\sigma,P_\sigma)\to0$ as
$v_\varphi\to\pm\infty$ and its $v_\varphi$ derivative is integrable.
Then
\begin{equation}
 c_{12}+c_{21}=0.
 \label{eq:offdiag-cancel}
\end{equation}
\end{lemma}

\begin{proof}
For each species,
\[
 \partial_{v_\varphi}\mu_K^\sigma(E_\sigma,P_\sigma)
 =\widehat v_\varphi\mu_{K,e}^\sigma(E_\sigma,P_\sigma)
  +r\mu_{K,p}^\sigma(E_\sigma,P_\sigma).
\]
Integration over $v_\varphi\in\R$ gives zero.  Summing the two identities
proves \eqref{eq:offdiag-cancel}.
\end{proof}

Assume now the physically standard energy monotonicity
\begin{equation}
 \mu_{K,e}^\sigma\le0,\qquad \sigma=\pm1.
 \label{eq:energy-monotone}
\end{equation}
Then
\[
 \mathcal E(K):=L_\phi-c_{22}
\]
is strictly positive with Dirichlet boundary data.  With $c=c_{12}$ and
$H=L_A-c_{11}$, Lemma \ref{lem:cancellation} writes
\begin{equation}
 \calL(K)=\begin{pmatrix}H&-c\\c&\mathcal E\end{pmatrix}.
 \label{eq:block-L}
\end{equation}

\begin{proposition}[Self-adjoint Schur-complement criterion]
\label{prop:Schur}
On the $L^2$ Hilbert realization define
\begin{equation}
 \calS(K)=L_A-c_{11}+c_{12}(L_\phi-c_{22})^{-1}c_{12}.
 \label{eq:Schur}
\end{equation}
Then $\calS(K)$ is self-adjoint with compact resolvent and
\begin{equation}
 \Ker\calL(K)=\left\{
 \binom{a}{-\mathcal E(K)^{-1}c_{12}a}:a\in\Ker\calS(K)
 \right\}.
 \label{eq:kernel-Schur}
\end{equation}
In particular, a necessary condition for bifurcation at $K_*$ is
$0\in\sigma(\calS(K_*))$.  If zero is a simple eigenvalue of $\calS(K_*)$,
then
\[
 \dim\Ker\calL(K_*)=\operatorname{codim}\Ran\calL(K_*)=1.
\]
If $a_*$ spans $\Ker\calS(K_*)$, kernel and adjoint-kernel vectors are
\begin{equation}
 \Psi=\binom{a_*}{-\mathcal E_*^{-1}c_*a_*},\qquad
 \Psi^*=\binom{a_*}{\mathcal E_*^{-1}c_*a_*}.
 \label{eq:kernel-vectors}
\end{equation}
\end{proposition}

\begin{proof}
The second row of $\calL(a,b)=0$ is
$ca+\mathcal Eb=0$, hence $b=-\mathcal E^{-1}ca$.  Substitution in the
first row gives $\calS a=0$.  Conversely this construction gives every
kernel vector.  The adjoint of \eqref{eq:block-L} is
$\left(\begin{smallmatrix}H&c\\-c&\mathcal E\end{smallmatrix}\right)$,
which gives the second vector in \eqref{eq:kernel-vectors}.  All
coefficient multipliers are real, $\mathcal E^{-1}$ is self-adjoint, and
$c\mathcal E^{-1}c$ is bounded and self-adjoint.  Thus $\calS$ is a
self-adjoint bounded perturbation of $L_A$ and has compact resolvent.

To relate this Hilbert realization to the strong Banach-space map, set
$P=\operatorname{diag}(L_A^{-1},L_\phi^{-1})$.  Then
\[
 D_z\calG(z_0(K),K)=P\calL(K).
\]
By \eqref{eq:compact}, this operator is identity minus compact and hence
Fredholm of index zero.  Elliptic bootstrapping identifies its kernel and
cokernel with those of the Hilbert realization above.  In particular,
the asserted cokernel dimension follows.  Elliptic regularity also puts
the displayed adjoint-kernel vector in $L^{p'}(\Omega)^2$; it therefore
defines the required bounded cokernel functional on the strong residual
space $\calY=L^p(\Omega)^2$.
\end{proof}


\subsection{Crossing form along a reference branch}

Suppose $K\mapsto z_0(K)$ is $C^1$.  A dot denotes total differentiation
along this branch.  For $\nu=e,p$,
\begin{equation}
 \dot\mu^\sigma_{K,\nu}
 =\partial_K\mu^\sigma_{K,\nu}
  +\sigma\dot\phi_0\mu^\sigma_{K,\nu e}
  +\sigma r\dot A_0\mu^\sigma_{K,\nu p}.
 \label{eq:profile-total-derivative}
\end{equation}
Replacing each $\mu^\sigma_{K,\nu}$ by
$\dot\mu^\sigma_{K,\nu}$ in \eqref{eq:c11c12}--\eqref{eq:c21c22}
defines $\dot c_{ij}$.  Differentiating \eqref{eq:Schur} yields
\begin{align}
 \dot\calS={}&-\dot c_{11}
 +\dot c_{12}\mathcal E^{-1}c_{12}
 +c_{12}\mathcal E^{-1}\dot c_{12}
 +c_{12}\mathcal E^{-1}\dot c_{22}\mathcal E^{-1}c_{12}.
 \label{eq:Sdot}
\end{align}
For the vectors in \eqref{eq:kernel-vectors},
\begin{equation}
 \left\langle\Psi^*,\dot\calL(K_*)\Psi\right\rangle
 =\left\langle a_*,\dot\calS(K_*)a_*\right\rangle.
 \label{eq:crossing-scalar}
\end{equation}
Thus the transversality condition is reduced to a scalar eigenvalue
crossing for a self-adjoint operator.

\subsection{Local secondary branches and the first branch coefficient}

We first state and prove the local principle.  The coefficient in
\eqref{eq:K1} gives the initial parameter direction; Section
\ref{sec:local-bif} reduces it to velocity moments.  Let $I\subset\R$ be open and suppose
\begin{equation}
 \calF(z_0(K),K)=0,\qquad K\in I,
 \label{eq:reference-branch}
\end{equation}
for a $C^2$ reference branch.  Introduce the shifted map
\[
 \widetilde{\calF}(x,K)=\calF(z_0(K)+x,K).
\]
Then $\widetilde{\calF}(0,K)=0$ for every $K$ and
$D_x\widetilde{\calF}(0,K)=\calL(K)$.

\begin{theorem}[Simple-eigenvalue local bifurcation]
\label{thm:CR}
Assume that $\calF$ is $C^3$ near $(z_0(K_*),K_*)$ and that
\begin{equation}
 \Ker\calL(K_*)=\operatorname{span}\{\Psi\},\qquad
 \Ker\calL(K_*)^*=\operatorname{span}\{\Psi^*\}.
 \label{eq:simple-kernel}
\end{equation}
Choose any nonzero adjoint-kernel vector $\Psi^*$ and suppose
\begin{equation}
 \left\langle\Psi^*,\dot\calL(K_*)\Psi\right\rangle\ne0.
 \label{eq:CR-transversality}
\end{equation}
Then $(z_0(K_*),K_*)$ is a bifurcation point.  In a neighborhood of this
point, all solutions consist of the reference branch and one additional
$C^2$ curve
\begin{align}
 z(s)&=z_0(K(s))+s\Psi+s^2\Xi+o(s^2),\label{eq:local-z}\\
 K(s)&=K_*+K_1s+o(s),\label{eq:local-K}
\end{align}
where $\Xi$ lies in a fixed complement of $\operatorname{span}\{\Psi\}$
and
\begin{equation}
 K_1=-\frac{
 \left\langle\Psi^*,D_z^2\calF(z_0(K_*),K_*)[\Psi,\Psi]\right\rangle}
 {2\left\langle\Psi^*,\dot\calL(K_*)\Psi\right\rangle}.
 \label{eq:K1}
\end{equation}
\end{theorem}

\begin{proof}[Proof of Theorem~\ref{thm:CR}]
The operator $\calL(K_*)$ is Fredholm of index zero by
\eqref{eq:compact}.  The first condition in \eqref{eq:simple-kernel}
therefore gives range codimension one.  Since
\[
 D_{xK}\widetilde{\calF}(0,K_*)\Psi
 =\dot\calL(K_*)\Psi,
\]
condition \eqref{eq:CR-transversality} is equivalent to
$D_{xK}\widetilde{\calF}(0,K_*)\Psi\notin\Ran\calL(K_*)$.
The Crandall--Rabinowitz theorem \cite{CrandallRabinowitz1971} applies.
Differentiating $\widetilde{\calF}(x(s),K(s))=0$ twice at $s=0$ and
pairing with $\Psi^*$ gives \eqref{eq:K1}.
\end{proof}

The numerator in \eqref{eq:K1} is completely explicit.  If
$h=(a,b)$ and $\widetilde h=(\widetilde a,\widetilde b)$, let
\begin{align}
 Q_\sigma(h,\widetilde h)
 ={}&\mu^\sigma_{K,ee}b\widetilde b
 +r\mu^\sigma_{K,ep}(b\widetilde a+a\widetilde b)
 +r^2\mu^\sigma_{K,pp}a\widetilde a.
 \label{eq:Qsigma}
\end{align}
Then
\begin{equation}
 D_z^2\calF[h,\widetilde h]
 =-\sum_{\sigma=\pm1}\sigma
 \binom{\displaystyle\int_{\R^3}\widehat v_\varphi
 Q_\sigma(h,\widetilde h)\dd v}
 {\displaystyle\int_{\R^3}Q_\sigma(h,\widetilde h)\dd v}.
 \label{eq:Fsecond}
\end{equation}
Equations \eqref{eq:Sdot}, \eqref{eq:crossing-scalar}, and
\eqref{eq:Fsecond} turn both the existence test and the leading branch
direction into velocity integrals and elliptic solves.

\subsection{The global bifurcation alternative}

Returning to the reference-branch setting of Theorem \ref{thm:CR}, let
$\mathscr U\subset\calX\times I$ be the open admissible domain of the
profiles and of $\calF$, containing the reference branch, and set
\[
 \mathscr S=\{(z,K)\in\mathscr U:\calF(z,K)=0\},\qquad
 \mathscr T=\{(z_0(K),K):K\in I\}.
\]
Let $\mathscr C_*$ be the connected component of
$\overline{\mathscr S\setminus\mathscr T}$ containing
$(z_0(K_*),K_*)$.

\begin{theorem}[Rabinowitz global alternative]
\label{thm:global-bif}
Assume the hypotheses of Theorem \ref{thm:CR} and that
\eqref{eq:compact} holds uniformly on bounded subsets of
$\calX\times I$.  Closures and boundaries below are taken in the ambient
space $\calX\times\R$.  Then at least one of the following alternatives
occurs:
\begin{enumerate}[label=\textup{(\alph*)},leftmargin=2.4em]
\item $\mathscr C_*$ is unbounded in $\calX\times I$;
\item $\mathscr C_*$ approaches $\partial\mathscr U$ (in particular, a
parameter or profile-admissibility boundary);
\item $\mathscr C_*$ meets the reference branch again at
$(z_0(K_1),K_1)$, where $K_1\ne K_*$ and $\calL(K_1)$ is singular.
\end{enumerate}
If $\mathscr U=\calX\times\R$, only
\textup{(a)} and \textup{(c)} can occur.
\end{theorem}

\begin{proof}[Proof of Theorem~\ref{thm:global-bif}]
Apply $\operatorname{diag}(L_A^{-1},L_\phi^{-1})$ and shift by the
reference branch.  The equation becomes
\[
 x-\mathcal K(x,K)=0,
\]
where $\mathcal K$ is compact on bounded sets and $x=0$ is a solution for
every $K$.  The transverse one-dimensional crossing in Theorem
\ref{thm:CR} has parity $-1$, equivalently odd crossing number, so the
local index changes sign.  Moreover, \eqref{eq:compact} makes every
bounded closed zero set that stays inside the admissible domain compact.
The parity/index-jump form of the Rabinowitz global theorem for
Fredholm maps of index zero
\cite{Rabinowitz1971,FitzpatrickPejsachowiczRabier1992,
DancerGlobal1973,BuffoniToland,Kielhoefer} therefore rules out a bounded
component that is separated from the admissibility boundary and from
every other singular point of the reference branch.  This gives the
three alternatives.
\end{proof}



\section{Computable local bifurcation criteria and explicit examples}
\label{sec:local-bif}

The abstract results of Section \ref{sec:linearization} involve the
kernel of \(\calL\) and an adjoint crossing form.  We reduce these
linearized objects to velocity moments and scalar elliptic eigenvalue
problems, obtaining tests that locate, count, and separate simple
crossings.  Branch direction additionally uses the nonlinear
coefficient in Theorem~\ref{thm:CR}, evaluated for particular profiles
in Section~\ref{sec:global-bif}.  

\subsection{Endpoint, monotonicity, and spectral-gap criteria}
\label{sec:explicit-bif-criteria}

We now replace the abstract requirements ``simple kernel'' and
``transversality'' by inequalities that can be checked from the velocity
moments and scalar elliptic eigenvalue problems.  Let
\[
 \mathcal V_A=\{a:e^{i\varphi}a\in H^1_0(\Omega;\mathbb C)\}
\]
be the magnetic form domain.  Let
$0<\lambda_{A,1}<\lambda_{A,2}\le\cdots$ be the eigenvalues of $L_A$ in
the real axisymmetric toroidal-amplitude space, and let $a_{A,1}$ be an
$L^2$-normalized first eigenfunction.  Thus
\begin{align}
 \lambda_{A,1}
 &=\inf_{a\in\mathcal V_A\setminus\{0\}}
 \frac{\displaystyle\int_\Omega
 (|\nabla a|^2+r^{-2}a^2)\dd x}{\|a\|_2^2},
 \label{eq:lambda-A}\\
 \lambda_E(K)
 &=\inf_{b\in H^1_0(\Omega)\setminus\{0\}}
 \frac{\displaystyle\int_\Omega
 (|\nabla b|^2-c_{22}(K)b^2)\dd x}{\|b\|_2^2}>0.
 \label{eq:lambda-E}
\end{align}
For the Schur complement \eqref{eq:Schur}, write
\begin{align}
 \mathfrak s_K[a]
 ={}&\int_\Omega
 \bigl(|\nabla a|^2+r^{-2}a^2-c_{11}(K)a^2\bigr)\dd x
 \notag\\
 &+\left\langle
 \mathcal E(K)^{-1}(c_{12}(K)a),c_{12}(K)a
 \right\rangle_{L^2},
 \label{eq:Schur-form}\\
 \nu_1(K)&=\inf_{a\ne0}\frac{\mathfrak s_K[a]}{\|a\|_2^2}.
 \label{eq:nu1}
\end{align}
Thus $\nu_1(K)$ is the lowest eigenvalue of $\calS(K)$.  Notice that the
last term in \eqref{eq:Schur-form} is nonnegative.

\begin{theorem}[Endpoint, monotonicity, and gap criterion]
\label{thm:computable-bif}
Let $z_0(K)$ be a $C^2$ reference branch on
$[K_-,K_+]$, let $\calF$ be $C^3$, and assume
\eqref{eq:energy-monotone}.  Suppose the coefficient functions are $C^1$
in $K$ with values in $L^\infty(\Omega)$, $\lambda_{A,1}$ is simple,
and the following conditions hold.
\begin{enumerate}[label=\textup{(C\arabic*)},leftmargin=2.8em]
\item For some $e_0>0$,
\begin{equation}
 \lambda_E(K)\ge e_0,
 \qquad K_-\le K\le K_+.
 \label{eq:C1-coercive}
\end{equation}
\item At the left endpoint,
\begin{equation}
 \|(c_{11}(K_-))_+\|_{L^\infty}<\lambda_{A,1},
 \label{eq:C2-left}
\end{equation}
where $g_+=\max\{g,0\}$.  At the right endpoint,
\begin{equation}
 \int_\Omega c_{11}(K_+)a_{A,1}^2\dd x
 >\lambda_{A,1}
 +\lambda_E(K_+)^{-1}
   \int_\Omega c_{12}(K_+)^2a_{A,1}^2\dd x.
 \label{eq:C2-right}
\end{equation}
\item There is $\kappa>0$ such that, for every
$K\in[K_-,K_+]$,
\begin{align}
 \operatorname*{ess\,inf}_{x\in\Omega}\dot c_{11}(K,x)
 \ge{}&\kappa
 +\frac{2\|\dot c_{12}(K)\|_\infty
              \|c_{12}(K)\|_\infty}{e_0}
 \notag\\
 &+\frac{\|(\dot c_{22}(K))_+\|_\infty
              \|c_{12}(K)\|_\infty^2}{e_0^2}.
 \label{eq:C3-monotone}
\end{align}
\item Uniformly on the interval,
\begin{equation}
 \sup_{K_-\le K\le K_+}
 \|(c_{11}(K))_+\|_\infty<\lambda_{A,2}.
 \label{eq:C4-gap}
\end{equation}
\end{enumerate}
Then there is a unique $K_*\in(K_-,K_+)$ for which
$\nu_1(K_*)=0$.  This eigenvalue is simple and
\begin{equation}
 \left\langle\Psi^*,\dot\calL(K_*)\Psi\right\rangle
 =\dot{\mathfrak s}_{K_*}[a_*]
 \le-\kappa\|a_*\|_2^2<0.
 \label{eq:explicit-crossing-negative}
\end{equation}
Consequently $(z_0(K_*),K_*)$ is a
Crandall--Rabinowitz bifurcation point, and the local conclusion of
Theorem \ref{thm:CR} applies.  Additional profile-specific hypotheses
that force field-norm escape are treated in
Section \ref{sec:global-bif}.
\end{theorem}

\begin{remark}
The local bifurcation condition can be viewed as a matching of a geometric mode length with a kinetic response length.  
\end{remark}

\begin{proof}
Condition \eqref{eq:C2-left}, the definition of $\lambda_{A,1}$, and
nonnegativity of the last term in \eqref{eq:Schur-form} give
\[
 \nu_1(K_-)\ge
 \lambda_{A,1}-\|(c_{11}(K_-))_+\|_\infty>0.
\]
Since $\|\mathcal E(K_+)^{-1}\|_{2\to2}=\lambda_E(K_+)^{-1}$,
\eqref{eq:C2-right} gives
$\mathfrak s_{K_+}[a_{A,1}]<0$, and hence
$\nu_1(K_+)<0$.  Form continuity and the min--max principle imply that
$\nu_1$ is continuous; see \cite{Kato}.  Thus it vanishes at least once.

For $a\in\mathcal V_A$, set
\[
 u_K[a]=\mathcal E(K)^{-1}(c_{12}(K)a).
\]
Formula \eqref{eq:Sdot} gives the exact form derivative
\begin{equation}
 \dot{\mathfrak s}_K[a]
 =-\int_\Omega\dot c_{11}a^2\dd x
 +2\int_\Omega\dot c_{12}a u_K[a]\dd x
 +\int_\Omega\dot c_{22}u_K[a]^2\dd x.
 \label{eq:form-derivative}
\end{equation}
By \eqref{eq:C1-coercive},
$\|u_K[a]\|_2\le e_0^{-1}\|c_{12}(K)\|_\infty\|a\|_2$.
Cauchy--Schwarz and \eqref{eq:C3-monotone} therefore yield
\begin{equation}
 \dot{\mathfrak s}_K[a]\le-\kappa\|a\|_2^2.
 \label{eq:form-strict-decrease}
\end{equation}
After integration in $K$, the min--max principle shows that $\nu_1$ is
strictly decreasing, so its zero is unique.

It remains to verify simplicity.  If $a\perp a_{A,1}$, then the
nonnegative Schur term and the min--max principle give
\[
 \mathfrak s_{K_*}[a]
 \ge\bigl(\lambda_{A,2}
 -\|(c_{11}(K_*))_+\|_\infty\bigr)\|a\|_2^2>0.
\]
Consequently the second eigenvalue of $\calS(K_*)$ is positive, whereas
its first eigenvalue is zero; hence the zero is simple.  For its
eigenfunction $a_*$, \eqref{eq:form-strict-decrease},
\eqref{eq:crossing-scalar}, and the Hellmann--Feynman formula give
\eqref{eq:explicit-crossing-negative}.  Proposition \ref{prop:Schur} and
Theorem \ref{thm:CR} finish the proof.
\end{proof}

\begin{proposition}[Counting and separation of crossings]
\label{prop:multi-crossing-count}
Let $[K_a,K_b]$ be a parameter interval on which
\eqref{eq:C1-coercive} holds, and write
\[
 \nu_1(K)\le\nu_2(K)\le\cdots
\]
for the eigenvalues of $\calS(K)$, repeated according to multiplicity.
Suppose we have bounds for the magnetic eigenvalues:
\begin{equation}
 \ell_j\le\lambda_{A,j}\le u_j,\qquad 1\le j\le J+1,
 \label{eq:multi-magnetic-enclosures}
\end{equation}
and define directly from the coefficients
\begin{align}
 m_{11}(K)&=\operatorname*{ess\,inf}_{x\in\Omega}c_{11}(K,x),
 &M_{11}(K)&=\operatorname*{ess\,sup}_{x\in\Omega}c_{11}(K,x),\notag\\
 q_{12}(K)&=\frac{\|c_{12}(K)\|_\infty^2}{e_0},
 &\underline\sigma_j(K)&=\ell_j-M_{11}(K),\notag\\
 &&\overline\sigma_j(K)&=u_j-m_{11}(K)+q_{12}(K).
 \label{eq:multi-coefficient-enclosures}
\end{align}
Suppose also that the following computable form bounds hold for every
$a\in\mathcal V_A$ and $K\in[K_a,K_b]$:
\begin{equation}
 -\kappa_+\|a\|_2^2
 \le \dot{\mathfrak s}_K[a]
 \le-\kappa_-\|a\|_2^2,
 \qquad 0<\kappa_-\le\kappa_+.
 \label{eq:multi-form-bounds}
\end{equation}
One coefficient-level choice is
\begin{align}
 \kappa_-={}&\inf_{K\in[K_a,K_b]}
 \left\{\operatorname*{ess\,inf}\dot c_{11}
 -\frac{2\|\dot c_{12}\|_\infty\|c_{12}\|_\infty}{e_0}
 -\frac{\|(\dot c_{22})_+\|_\infty\|c_{12}\|_\infty^2}{e_0^2}
 \right\}>0,                                      \label{eq:multi-kappa-minus}\\
 \kappa_+={}&\sup_{K\in[K_a,K_b]}
 \left\{\|\dot c_{11}\|_\infty
 +\frac{2\|\dot c_{12}\|_\infty\|c_{12}\|_\infty}{e_0}
 +\frac{\|\dot c_{22}\|_\infty\|c_{12}\|_\infty^2}{e_0^2}
 \right\}.                                        \label{eq:multi-kappa-plus}
\end{align}
Then the coefficient quantities in \eqref{eq:multi-coefficient-enclosures}
automatically enclose the Schur eigenvalues:
\begin{equation}
 \underline\sigma_j(K)\le\nu_j(K)\le\overline\sigma_j(K),
 \qquad 1\le j\le J+1.                         \label{eq:multi-derived-enclosures}
\end{equation}
Consequently every index $1\le j\le J$ satisfying
\begin{equation}
 \underline\sigma_j(K_a)>0,\qquad
 \overline\sigma_j(K_b)<0                     \label{eq:multi-certified-index}
\end{equation}
has exactly one crossing parameter $K_j\in(K_a,K_b)$.  It belongs to
the computable interval $I_j=[L_j,U_j]$, where
\begin{align}
 L_j&=\max\left\{
 K_a+\frac{\underline\sigma_j(K_a)}{\kappa_+},
 K_b+\frac{\underline\sigma_j(K_b)}{\kappa_-}\right\},\notag\\
 U_j&=\min\left\{
 K_a+\frac{\overline\sigma_j(K_a)}{\kappa_-},
 K_b+\frac{\overline\sigma_j(K_b)}{\kappa_+}\right\}.
 \label{eq:multi-crossing-enclosure}
\end{align}
For consecutive certified indices,
\begin{equation}
 \max\{0,L_{j+1}-U_j\}
 \le K_{j+1}-K_j
 \le U_{j+1}-L_j.                                \label{eq:multi-gap-enclosure}
\end{equation}
In particular, $L_{j+1}>U_j$ certifies a strictly positive gap.

The number of crossing eigenvalues, counted with multiplicity, admits
the finite bounds
\begin{equation}
 \#\{1\le j\le J:\underline\sigma_j(K_a)>0,\
                      \overline\sigma_j(K_b)<0\}
 \le N_{\rm cross}\le
 \#\{1\le j\le J:\overline\sigma_j(K_a)>0,\
                      \underline\sigma_j(K_b)<0\},
 \label{eq:multi-count-bounds}
\end{equation}
provided the coefficient bounds include every mode whose derived
enclosure can meet zero.  A convenient tail certificate is
$\underline\sigma_1(K_a)>0$ together with
$\underline\sigma_{J+1}(K_b)>0$: then no index above $J$ can cross.  If,
in addition, each interval
$[\underline\sigma_j(K_b),\overline\sigma_j(K_b)]$, $1\le j\le J$,
avoids zero, the exact count is
\begin{equation}
 N_{\rm cross}
 =\#\{1\le j\le J:\overline\sigma_j(K_b)<0\}.
                                                        \label{eq:multi-exact-count}
\end{equation}
For each adjacent index $k=j\pm1$, suppose that the same coefficient
data either exclude a zero of $\nu_k$ on $[K_a,K_b]$ or confine its
unique zero to an interval $I_k$ disjoint from $I_j$ (with the evident
one-sided convention for $j=1$).  Then zero is a simple eigenvalue of
$\calS(K_j)$; hence $K_j$ is a Crandall--Rabinowitz bifurcation point.
Thus \eqref{eq:multi-gap-enclosure} estimates gaps between genuine
local bifurcation points whenever simplicity is certified.
\end{proposition}

\begin{proof}
Coercivity of $\mathcal E(K)$ and the Schur form \eqref{eq:Schur-form}
give, for every $a\in\mathcal V_A$,
\begin{align*}
 \langle L_Aa,a\rangle-M_{11}(K)\|a\|_2^2
 &\le\mathfrak s_K[a]\\
 &\le\langle L_Aa,a\rangle+
 \bigl[-m_{11}(K)+q_{12}(K)\bigr]\|a\|_2^2.
\end{align*}
The min--max principle and \eqref{eq:multi-magnetic-enclosures} prove
\eqref{eq:multi-derived-enclosures}; this is the step that removes any
Schur-eigenvalue input from the hypotheses.

Integrating \eqref{eq:multi-form-bounds} and applying the min--max
principle gives, for $K_a\le K_1<K_2\le K_b$,
\begin{equation}
 -\kappa_+(K_2-K_1)
 \le\nu_j(K_2)-\nu_j(K_1)
 \le-\kappa_-(K_2-K_1).                          \label{eq:multi-eigenvalue-slope}
\end{equation}
Thus every ordered eigenvalue is strictly decreasing and can vanish at
most once.  The derived endpoint signs give existence.  Applying
\eqref{eq:multi-eigenvalue-slope} from each endpoint yields
\eqref{eq:multi-crossing-enclosure}; subtraction gives
\eqref{eq:multi-gap-enclosure}.  Monotonicity also shows that the number
of crossings is the increase of the Morse index from $K_a$ to $K_b$,
which proves \eqref{eq:multi-count-bounds} and
\eqref{eq:multi-exact-count}.  The tail certificate follows from the
ordering of the magnetic and Schur eigenvalues.  Excluding adjacent
zeros, or confining them to disjoint intervals, excludes multiplicity.
The form derivative is then nonzero on the kernel by
\eqref{eq:multi-form-bounds}, so Proposition \ref{prop:Schur} and
Theorem \ref{thm:CR} apply.  Finally,
\eqref{eq:multi-kappa-minus}--\eqref{eq:multi-kappa-plus} follow from
\eqref{eq:form-derivative}, coercivity, and Cauchy--Schwarz.
\end{proof}

\subsection{General profiles, torus geometry, and multiple crossings}
\label{sec:explicit-C1-C4}

We now verify the local crossing conditions directly for broad
two-species profile classes.  The presentation begins with the general
even--odd decomposition, then supplies computable variable-torus
spectral bounds and a catalogue of admissible energy--momentum laws.
This subsumes the product examples and avoids repeating a special
profile as a separate benchmark.  The multiple-crossing criteria
from Proposition \ref{prop:multi-crossing-count} will be evaluated for
the same profiles and geometries.

\subsubsection*{General even--odd profiles}

The construction does not require the energy
law $e^{-e}$, a product ansatz, or the amplitude $K^2$.  We first separate
the profile assumptions from the geometric spectral window.
Throughout this part, $\Omega$ is connected and satisfies
\eqref{eq:domain}, and
$0<r_-\le r\le r_+<\infty$ on $\overline\Omega$.
Let $M_{\rm ev}(e,p)$ be even in $p$ and $M_{\rm od}(e,p)$ odd in $p$;
these subscripts refer to parity, not to differentiation in energy.
Assume they have the analytic moment bounds of Assumption
\ref{ass:analytic}, separately for each profile, and that for some $B\ge1$,
\begin{align}
 M_{\rm ev}&>0,\qquad |M_{\rm od}|\le B M_{\rm ev},
       &&(e,p)\in\R^2,\notag\\
 -\partial_eM_{\rm ev}&>0,\qquad
 |\partial_eM_{\rm od}|\le B(-\partial_eM_{\rm ev}),
       &&e\ge\tfrac12,\ p\in\R.                 \label{eq:general-profile-dom}
\end{align}
The threshold $1/2$ only supplies a neighborhood of the zero-field
branch: there $e_\pm=\av v\ge1$.  Global energy monotonicity for negative
energies is not needed for the local crossing argument.
Define continuous velocity moments on $[r_-,r_+]$ by
\begin{align}
 W(r)&=2r\int_{\R^3}\widehat v_\varphi
       \partial_pM_{\rm ev}(\av v,rv_\varphi)\dd v,
                                                        \label{eq:general-W}\\
 Q(r)&=-2\int_{\R^3}\widehat v_\varphi
       \partial_eM_{\rm od}(\av v,rv_\varphi)\dd v,
                                                        \label{eq:general-Q}\\
 D(r)&=-2\int_{\R^3}
       \partial_eM_{\rm ev}(\av v,rv_\varphi)\dd v>0.
                                                        \label{eq:general-D}
\end{align}
We require $W(r)>0$ on the radial interval, and choose bounds
\[
 0<W_-\le W(r)\le W_+,\qquad Q_\infty=\max_{[r_-,r_+]}|Q(r)|.
\]
Only these three moments enter the linear crossing criterion.

\begin{theorem}[A general class satisfying \textup{(C1)--(C4)}]
\label{thm:general-profile-class}
Assume \eqref{eq:general-profile-dom} and the moment conditions above.
Let $e_0>0$ be any lower bound for the first Dirichlet eigenvalue of
$L_\phi$, and choose spectral bounds
\begin{equation}
 0<\ell_1\le\lambda_{A,1}\le u_1,\qquad
 \lambda_{A,2}\ge\ell_2,\qquad
 \frac{W_+}{W_-}<\frac{\ell_2}{u_1}.             \label{eq:general-window}
\end{equation}
Set
\begin{align}
 t_+&=\frac12\left(\frac{u_1}{W_-}+\frac{\ell_2}{W_+}\right),
 &t_-&=\frac{\ell_1}{2W_+},
 &\Delta&=t_+W_--u_1>0.                         \label{eq:general-endpoints}
\end{align}
Let $\eta$ be real analytic on an open interval containing
$[K_-,K_+]$, with
\begin{equation}
 \eta(K_\pm)=t_\pm,\qquad
 \eta'(K)\ge\eta_0>0\quad(K_-\le K\le K_+).
                                                        \label{eq:general-amplitude}
\end{equation}
On this interval define
\begin{equation}
 \mu_{K,\epsilon}^+(e,p)=\mu_{K,\epsilon}^-(e,p)
   =\eta(K)\bigl(M_{\rm ev}(e,p)+\epsilon M_{\rm od}(e,p)\bigr).
                                                        \label{eq:general-profiles}
\end{equation}
For
\begin{equation}
 |\epsilon|<\epsilon_{\rm gen}:=
 \min\left\{\frac1{2B},
       \sqrt{\frac{e_0\Delta}{2t_+^2Q_\infty^2}},
       \sqrt{\frac{e_0W_-}{4t_+Q_\infty^2}}\right\},       \label{eq:general-eps}
\end{equation}
where terms with $Q_\infty=0$ are omitted, both species are positive,
$(A,\phi)=(0,0)$ is an exact reference equilibrium for all $K$, and
\textup{(C1)--(C4)} hold with
\begin{equation}
 \kappa=\frac12\eta_0W_-.
                                                        \label{eq:general-kappa}
\end{equation}
In particular, there is a unique simple transverse stationary crossing
$K_*\in(K_-,K_+)$.  The local bifurcation theorem applies.  Section
\ref{sec:global-bif} gives separate norm-escape mechanisms for
specific two-species profile families.

For $\epsilon=0$, the crossing is
$\eta(K_*)=\Lambda_1(W)$ and its kernel is purely magnetic.
For $\epsilon\ne0$ and $Q(r)>0$ on $[r_-,r_+]$, it is genuinely coupled:
\begin{equation}
 \Psi=\binom{a_*}
 {\epsilon\eta(K_*)(L_\phi+\eta(K_*)D)^{-1}(Qa_*)},
 \qquad a_*\ne0,                               \label{eq:general-kernel}
\end{equation}
and both components are nonzero.
\end{theorem}

\begin{proof}
The pointwise bounds \eqref{eq:general-profile-dom} imply positivity and
energy decrease near the zero branch.  Equal species cancel both charge
and current at $A=\phi=0$, whether or not their common profile is even.
Parity and Lemma \ref{lem:cancellation} give the exact formulas
\begin{equation}
 c_{11}=\eta W,\qquad c_{12}=-\epsilon\eta Q,\qquad
 c_{21}=\epsilon\eta Q,\qquad c_{22}=-\eta D.     \label{eq:general-coefficients}
\end{equation}
For completeness, the odd part contributes nothing to $c_{11}$ since
$\partial_pM_{\rm od}$ is even; the even part contributes nothing to
$c_{12}$ since $\partial_eM_{\rm ev}$ is even.
The odd part integrates to zero in $c_{22}$.
The identity
$2r\int\partial_pM_{\rm od}\dd v
 =-2\int\widehat v_\varphi\partial_eM_{\rm od}\dd v=Q$
also directly verifies the sign of $c_{21}$.

Thus $\mathcal E=L_\phi+\eta D\ge L_\phi\ge e_0$, giving
\textup{(C1)}.  At the left endpoint,
$t_-W_+=\ell_1/2<\lambda_{A,1}$.  At the right endpoint,
for the normalized first magnetic eigenfunction,
\begin{align*}
 \int_\Omega c_{11}(K_+)a_{A,1}^2\dd x
 &\ge t_+W_-=u_1+\Delta,\\
 \lambda_E(K_+)^{-1}\int_\Omega c_{12}(K_+)^2a_{A,1}^2\dd x
 &\le\epsilon^2t_+^2Q_\infty^2/e_0<\Delta/2.
\end{align*}
These are \textup{(C2)}.  The upper window endpoint gives
$t_+W_+<\ell_2\le\lambda_{A,2}$, hence \textup{(C4)}.

To check \textup{(C3)}, keep the common factor $\eta'(K)$ rather than
bounding its maximum independently:
\begin{align}
 \dot c_{11}&=\eta'W,\qquad
 \dot c_{12}=-\epsilon\eta'Q,\qquad
 \dot c_{22}=-\eta'D\le0,\notag\\
 \frac{2\|\dot c_{12}\|_\infty\|c_{12}\|_\infty}{e_0}
 &\le\eta'\frac{2\epsilon^2t_+Q_\infty^2}{e_0}
 <\tfrac12\eta'W_-.                            \label{eq:general-C3}
\end{align}
Consequently the margin is at least $\eta_0W_-/2$.
The first magnetic eigenvalue is simple: on a connected meridional
domain away from the axis the operator is unitarily equivalent to
$-\Delta_{r,z}+3/(4r^2)$ and has a positive principal eigenfunction
\cite{Amann,GilbargTrudinger}.  Theorem \ref{thm:computable-bif} now applies.
The kernel formula is Proposition \ref{prop:Schur}; if $Q>0$, multiplication
by $Q$ cannot annihilate $a_*$, and the electric inverse is injective.
\end{proof}

\begin{proposition}[Mode counts for the general even--odd profiles]
\label{prop:general-profile-count}
Retain the profile and moment assumptions above, but let
$[K_a,K_b]$ be any interval on which
\begin{equation}
 0\le\eta_a:=\eta(K_a)<\eta(K_b)=:\eta_b\le T,
 \qquad
 0<\eta'_-\le\eta'(K)\le\eta'_+ .                \label{eq:general-count-amplitude}
\end{equation}
Put $D_+=\max_{[r_-,r_+]}D(r)$ and assume
\begin{equation}
 |\epsilon|<\frac1{2B},\qquad
 m_\epsilon:=
 W_--\frac{2\epsilon^2TQ_\infty^2}{e_0}>0.       \label{eq:general-count-smallness}
\end{equation}
If $\ell_j\le\lambda_{A,j}\le u_j$, then every Schur eigenvalue has the
computable enclosure
\begin{equation}
 \ell_j-\eta(K)W_+
 \le\nu_j(K)\le
 u_j-\eta(K)W_-+
 \frac{\epsilon^2\eta(K)^2Q_\infty^2}{e_0}.      \label{eq:general-count-enclosure}
\end{equation}
Moreover, \eqref{eq:multi-form-bounds} holds with
\begin{align}
 \kappa_-&=\eta'_-m_\epsilon,\notag\\
 \kappa_+&=\eta'_+\left(
 W_++\frac{2\epsilon^2TQ_\infty^2}{e_0}
 +\frac{\epsilon^2T^2D_+Q_\infty^2}{e_0^2}
 \right).                                       \label{eq:general-count-kappas}
\end{align}
Consequently Proposition \ref{prop:multi-crossing-count} applies using
the endpoint values in \eqref{eq:general-count-enclosure}.  In
particular, every index in
\begin{equation}
 \mathcal J_{\rm cert}:=
 \left\{j:
 \ell_j-\eta_aW_+>0,\quad
 u_j-\eta_bW_-+
 \frac{\epsilon^2\eta_b^2Q_\infty^2}{e_0}<0
 \right\}                                       \label{eq:general-certified-modes}
\end{equation}
has one and only one crossing in $(K_a,K_b)$, and its crossing interval
and the gaps to neighboring crossings are given by
\eqref{eq:multi-crossing-enclosure}--\eqref{eq:multi-gap-enclosure}.
Thus the count and gap calculation requires only the velocity moments
$W_\pm,Q_\infty,D_+$, the electric bound $e_0$, the scalar amplitude
bounds, and finitely many magnetic eigenvalue enclosures.

For $\epsilon=0$, the Schur operator is $L_A-\eta(K)W$, so its zeros
obey the exact weighted-mode relation $\eta(K_j)=\Lambda_j(W)$, where
$L_Aa=\Lambda_j(W)Wa$ and the eigenvalues are counted with multiplicity.
For $\epsilon\ne0$ and $Q(r)>0$, they give a genuinely coupled two-field count
without first locating the individual zeros numerically.
\end{proposition}

\begin{proof}
The coefficient identities \eqref{eq:general-coefficients} give
\[
 \calS(K)=L_A-\eta(K)W+
 \epsilon^2\eta(K)^2
 Q(L_\phi+\eta(K)D)^{-1}Q .
\]
The last term is nonnegative and has operator norm at most
$\epsilon^2\eta(K)^2Q_\infty^2/e_0$.  The min--max principle gives
\eqref{eq:general-count-enclosure}.  Formula
\eqref{eq:form-derivative}, together with
$\|(L_\phi+\eta D)^{-1}\|\le e_0^{-1}$, yields
\[
 \dot{\mathfrak s}_K[a]
 \le-\eta'(K)\left(
 W_--\frac{2\epsilon^2TQ_\infty^2}{e_0}\right)\|a\|_2^2.
\]
Bounding all three terms in absolute value gives the second constant in
\eqref{eq:general-count-kappas}.  The result follows from Proposition
\ref{prop:multi-crossing-count}.
\end{proof}

\begin{example}[A quadratic-amplitude coefficient table]
\label{ex:quadratic-coefficient-table}
Take the general even--odd profiles with
$\eta(K)=cK^2$, $c>0$, on $0<K_a<K_b$, and set
$T=cK_b^2$.  Proposition \ref{prop:multi-crossing-count} requires no
spectral calculation for $\calS(K)$.  Its complete endpoint input is
the following table, obtained from velocity moments and the fixed
magnetic eigenvalue bounds:
\begin{equation}
\begin{array}{c|cc}
 &\underline\sigma_j(K)&\overline\sigma_j(K)\\ \hline
 K=K_a&\ell_j-cK_a^2W_+
 &u_j-cK_a^2W_-+\epsilon^2c^2K_a^4Q_\infty^2/e_0\\[1mm]
 K=K_b&\ell_j-cK_b^2W_+
 &u_j-cK_b^2W_-+\epsilon^2c^2K_b^4Q_\infty^2/e_0.
\end{array}                                      \label{eq:quadratic-endpoint-table}
\end{equation}
If
\begin{equation}
 m_\epsilon=W_--\frac{2\epsilon^2cK_b^2Q_\infty^2}{e_0}>0,
 \label{eq:quadratic-meps}
\end{equation}
then admissible form-slope constants are
\begin{align}
 \kappa_-&=2cK_a m_\epsilon,\notag\\
 \kappa_+&=2cK_b\left(
 W_++\frac{2\epsilon^2cK_b^2Q_\infty^2}{e_0}
 +\frac{\epsilon^2c^2K_b^4D_+Q_\infty^2}{e_0^2}
 \right).                                      \label{eq:quadratic-kappas}
\end{align}
Thus the signs in \eqref{eq:quadratic-endpoint-table} count the
crossings, while substituting the same four entries and
\eqref{eq:quadratic-kappas} into
\eqref{eq:multi-crossing-enclosure}--\eqref{eq:multi-gap-enclosure}
encloses their locations and mutual gaps.  Every input is a velocity
quadrature, an elementary coefficient extremum, or an eigenvalue bound
for the parameter-independent operator $L_A$.
\end{example}

The amplitude requirement is an actual restriction: $\eta$ must reach
$t_+$ with a positive lower bound for $\eta'$ on the chosen interval.
For example, $\eta(K)=cK^{2n}$, $c>0$, $n\in\N$, works with
$K_\pm=(t_\pm/c)^{1/(2n)}$.  So does $\eta(K)=c(e^K-1)^2$ on $K>0$,
with $K_\pm=\log(1+\sqrt{t_\pm/c})$.
These examples also have $\eta(0)=\eta'(0)=0$, as in the original
parameter convention.  A bounded amplitude whose range misses the
spectral window is not covered.

\subsubsection*{A variable-torus criterion and mode-count example}

\begin{proposition}[Spectral, counting, and gap bounds on variable tori]
\label{prop:general-torus}
Let $R>0$, $0<a<R$, and
\begin{equation}
 \Omega_{R,a}
 =\{(r,\varphi,z):(r-R)^2+z^2<a^2\}.             \label{eq:variable-torus-domain}
\end{equation}
The spectral enclosures apply in the full axisymmetric space or in a
fixed $z$-reflection class.  A nonlinear bifurcation assertion in a
restricted class also requires that the full residual preserve that
class.  For the profiles below the domain and coefficients are
$z$-reflection invariant, so the $z$-even class is invariant; an
arbitrary $z$-odd restriction need not be invariant without an
additional oddness property of the nonlinearity.  Let
$0<\gamma_1\le\gamma_2\le\cdots$ be the corresponding Dirichlet
eigenvalues of the unit meridional disk, repeated with multiplicity.
In the full space these are the numbers $j_{m,n}^2$ with their usual
disk multiplicities.

Assume $W$ is positive and $C^1$ on $[R-a,R+a]$ and, for some $L\ge0$,
\begin{equation}
 |rW'(r)|\le L W(r).                              \label{eq:general-logW}
\end{equation}
With $s=a/R$ and $\zeta=-\log(1-s)$, set
\begin{equation}
 W_-=W(R)e^{-L\zeta},\qquad W_+=W(R)e^{L\zeta}.   \label{eq:variable-torus-W}
\end{equation}
For every $j\ge1$, valid magnetic eigenvalue enclosures are
\begin{equation}
 \ell_j(R,a):=\frac{\gamma_j}{a^2}
             +\frac{3}{4(R+a)^2}
 \le\lambda_{A,j}\le
 \frac{\gamma_j}{a^2}
             +\frac{3}{4(R-a)^2}
 =:u_j(R,a),                                     \label{eq:variable-torus-mode-enclosures}
\end{equation}
and a valid electric coercivity bound is
\begin{equation}
 e_0(R,a)=\frac{j_{0,1}^2}{a^2}
          -\frac{1}{4(R-a)^2}.                   \label{eq:variable-torus-e0}
\end{equation}
Thus $e_0>0$ is an explicit aspect-ratio test.

In the decoupled general even--odd class with
$\eta(K)=cK^{2n}$, $c>0$, the $j$th weighted crossing obeys
\begin{equation}
 \left(\frac{\ell_j(R,a)}{cW_+}\right)^{1/(2n)}
 \le K_j\le
 \left(\frac{u_j(R,a)}{cW_-}\right)^{1/(2n)}.     \label{eq:variable-torus-Kj}
\end{equation}
The following elementary inequality certifies a positive gap between
the $j$th and $(j+1)$st crossing enclosures:
\begin{equation}
 (1-s)^{-2L}<
 \frac{\gamma_{j+1}/s^2+3/[4(1+s)^2]}
      {\gamma_j/s^2+3/[4(1-s)^2]}.               \label{eq:variable-torus-gap-test}
\end{equation}
When it holds,
\begin{align}
 K_{j+1}-K_j\ge c^{-1/(2n)}
 \bigg[
 &\left(\frac{\ell_{j+1}(R,a)}{W_+}\right)^{1/(2n)}
 -\left(\frac{u_j(R,a)}{W_-}\right)^{1/(2n)}
 \bigg]>0.                                       \label{eq:variable-torus-gap}
\end{align}
Finally, the number of crossings in $0<K<K_{\max}$ satisfies
\begin{align}
 \#\left\{j:\frac{u_j(R,a)}{W_-}<cK_{\max}^{2n}\right\}
 \le N(K_{\max})
 \le
 \#\left\{j:\frac{\ell_j(R,a)}{W_+}<cK_{\max}^{2n}\right\}.
                                                        \label{eq:variable-torus-count}
\end{align}
If no enclosure straddles $cK_{\max}^{2n}$, the two counts agree.
Consecutive gap tests certify that the corresponding crossings are
simple Crandall--Rabinowitz bifurcation points.

For the coupled general even--odd class, use the same
$\ell_j,u_j,W_\pm,e_0$ in Proposition
\ref{prop:general-profile-count}; this gives the analogous count and
parameter-gap enclosures with the explicit $\epsilon^2$ correction.
\end{proposition}

\begin{proof}
In the axisymmetric $L^2$ measure $2\pi r\dd r\dd z$, the unitary
substitution $a\mapsto r^{1/2}a$ conjugates $L_A$ to
$-\Delta_{r,z}+3/(4r^2)$ on the meridional disk, while the electric
operator is conjugate to $-\Delta_{r,z}-1/(4r^2)$.  The disk spectrum,
the bounds $R-a\le r\le R+a$, and the min--max principle give
\eqref{eq:variable-torus-mode-enclosures} and
\eqref{eq:variable-torus-e0}.  Integrating
\eqref{eq:general-logW} gives \eqref{eq:variable-torus-W}.
The weighted min--max principle gives
$\lambda_{A,j}/W_+\le\Lambda_j(W)\le\lambda_{A,j}/W_-$;
inserting the geometric bounds and solving
$cK_j^{2n}=\Lambda_j(W)$ gives
\eqref{eq:variable-torus-Kj}--\eqref{eq:variable-torus-count};
\eqref{eq:variable-torus-gap-test} is precisely
$u_j/W_-<\ell_{j+1}/W_+$.
\end{proof}

\begin{example}[General even--odd profiles on a variable torus]
\label{ex:variable-torus-count}
In the catalogue below take
\[
 H(e)=e^{-\beta e},\qquad
 g(p)=1+\vartheta\bigl[1-(1+\omega^2p^2)^{-\nu}\bigr],
 \qquad \eta(K)=cK^2,
\]
with $\beta,\vartheta,\omega,\nu,c>0$.  First set $h=0$, so
$\epsilon=0$.  Formula \eqref{eq:catalogue-W} gives $W(R)$ by a single
positive velocity quadrature, while
$L=2\max\{1,\nu\}$ gives $W_\pm$ from
\eqref{eq:variable-torus-W}.  In the $z$-even reflection class the first
three disk numbers are
\[
 \gamma_1=j_{0,1}^2,\qquad
 \gamma_2=j_{1,1}^2,\qquad
 \gamma_3=j_{2,1}^2.
\]
Thus \eqref{eq:variable-torus-Kj} gives three explicit crossing
intervals, \eqref{eq:variable-torus-gap-test} decides whether the first
two gaps are certified, and \eqref{eq:variable-torus-count} counts all
crossings below a prescribed $K_{\max}$.  Every quantity is a Bessel
zero, a velocity integral, or one of the two geometric parameters
$(R,a)$.

For a numerical instance, take
\begin{equation}
 R=1,\qquad a=\frac1{10},\qquad
 \beta=\vartheta=\omega=\nu=c=1.                \label{eq:variable-torus-numerical-data}
\end{equation}
Then $L=2$, $W_-=0.81W(1)$, $W_+=0.9^{-2}W(1)$, and
$e_0>578$.  For reproducibility, the needed velocity moment reduces to
one dimension:
\begin{equation}
 W(1)=16\pi\int_0^\infty
 \frac{t^2e^{-\sqrt{1+t^2}}}{(1+t^2)^2}\dd t,
 \quad 0\le 16\pi\int_T^\infty
 \frac{t^2e^{-\sqrt{1+t^2}}}{(1+t^2)^2}\dd t
 \le\frac{16\pi e^{-T}}{T^2}\quad(T\ge1).
 \label{eq:variable-torus-moment-tail}
\end{equation}
Interval quadrature on $[0,T]$ and this explicit tail bound certify
any desired enclosure for $W(1)$.  The numerical intervals below
use units of $W(1)^{-1/2}$, so their endpoints require only Bessel
zeros, not a decimal value of $W(1)$.  For example, the standard
zero indexing and alternating-series sign checks for $J_m$ give
\[
\begin{array}{c|c}
 j_{0,1}&[2.4048255576,2.4048255578]\\
 j_{1,1}&[3.8317059701,3.8317059703]\\
 j_{2,1}&[5.1356223017,5.1356223020].
\end{array}
\]
Using these outward enclosures \cite{NISTDLMF} in
\eqref{eq:variable-torus-Kj}, the first three derived crossing
intervals in the $z$-even class, expressed in units of
$W(1)^{-1/2}$, are
\begin{equation}
 I_1\subset[21.655,26.742],\qquad
 I_2\subset[34.492,42.588],\qquad
 I_3\subset[46.226,57.073].                    \label{eq:variable-torus-numerical-intervals}
\end{equation}
Both instances of \eqref{eq:variable-torus-gap-test} hold: its left
side is $0.9^{-4}=1.5242$, while the first two right sides are $2.5357$
and $1.7957$.  Hence
\begin{equation}
 K_2-K_1\ge7.750\,W(1)^{-1/2},\qquad
 K_3-K_2\ge3.638\,W(1)^{-1/2}.                 \label{eq:variable-torus-numerical-gaps}
\end{equation}
Moreover, choosing $K_{\max}=45W(1)^{-1/2}$ makes the two counts in
\eqref{eq:variable-torus-count} equal to two.  Thus this single finite
calculation certifies exactly two simple bifurcation points below
$K_{\max}$ and a third one above it \emph{in the invariant $z$-even
class}.  In the unrestricted axisymmetric space some disk modes are
repeated, so these restricted simple-eigenvalue conclusions do not
assert unrestricted simplicity.

Here is a genuinely coupled application of the updated Proposition
\ref{prop:multi-crossing-count}.  Rescale the common profile so that
$W(1)=1$, retain the numerical data
\eqref{eq:variable-torus-numerical-data}, choose any nonzero odd factor
$h$ from \eqref{eq:catalogue-odd}, and take a nonzero $\epsilon$ so small
that
\begin{equation}
 |\epsilon|<\frac12,\qquad
 \epsilon^2Q_\infty^2\le10^{-6},\qquad
 \epsilon^2D_+Q_\infty^2\le10^{-4}.            \label{eq:variable-coupled-smallness}
\end{equation}
This is always possible because the two velocity moments are finite,
and $Q>0$ makes the resulting kernel genuinely two-field.  Apply the
quadratic table of Example \ref{ex:quadratic-coefficient-table} on
$[K_a,K_b]=[20,45]$.  Using $e_0=578$ and rounding outward, the derived
coefficient enclosures are
\begin{equation}
\begin{array}{c|cc}
 j&[\underline\sigma_j(20),\overline\sigma_j(20)]
  &[\underline\sigma_j(45),\overline\sigma_j(45)]\\ \hline
 1&[85.10,255.26]&[-1921.08,-1060.99]\\
 2&[974.98,1145.14]&[-1031.20,-171.11]\\
 3&[2144.24,2314.41]&[138.07,998.16].
\end{array}                                      \label{eq:variable-coupled-table}
\end{equation}
The same data allow the choices
\begin{equation}
 \kappa_-=32.39,\qquad \kappa_+=111.23.         \label{eq:variable-coupled-kappas}
\end{equation}
Indeed, these follow immediately from \eqref{eq:quadratic-kappas},
$W_-=0.81$, $W_+=0.9^{-2}$, and
\eqref{eq:variable-coupled-smallness}.  Therefore Proposition
\ref{prop:multi-crossing-count} gives exactly two crossings in
$(20,45)$ in the invariant $z$-even class and the fully
coefficient-level enclosures
\begin{equation}
 I_1\subset[20.765,27.881],\qquad
 I_2\subset[28.765,43.462],\qquad
 K_2-K_1>0.88.                                  \label{eq:variable-coupled-crossings}
\end{equation}
The intervals are disjoint, while the third mode remains positive at
$K=45$; hence both crossings are simple Crandall--Rabinowitz points
within that class.
No eigenvalue of $\calS(K)$ has been used as input.  For each fixed
finite list of distinct disk modes, analogous tests hold for all
sufficiently small $s=a/R$; repeated disk eigenvalues require either a
reflection-class restriction or sharper validated magnetic
eigenvalue enclosures, exactly as the simplicity criterion predicts.
\end{example}

Even without \eqref{eq:general-logW}, continuity and $W(R)>0$ imply
$\max_{[R-a,R+a]}W/\min_{[R-a,R+a]}W\to1$ as $a\downarrow0$.
Thus sufficiently thin variable tori work for every fixed profile with
such a positive continuous moment.  The logarithmic bound turns this
qualitative fact into the finite tests
\eqref{eq:variable-torus-gap-test} and
\eqref{eq:variable-torus-count}.

\subsubsection*{Explicit energy laws and nonseparable mixtures}

\begin{proposition}[A catalogue of admissible profiles]
\label{prop:profile-catalogue}
Let $J\in\N$ and let $a_j,\vartheta_j,\omega_j,\nu_j$ be positive parameters.
Choose
\begin{align}
 M_{\rm ev}(e,p)&=\sum_{j=1}^J a_jH_j(e)g_j(p),
 &M_{\rm od}(e,p)&=\sum_{j=1}^J a_jH_j(e)h_j(p), \label{eq:catalogue-mixture}\\
 g_j(p)&=1+\vartheta_j\left[1-(1+\omega_j^2p^2)^{-\nu_j}\right].
                                                        \label{eq:catalogue-even}
\end{align}
Independently for each $j$, take either energy law
\begin{equation}
 H_j(e)=e^{-\beta_je},\quad\beta_j>0,
 \qquad\hbox{or}\qquad
 H_j(e)=(b_j+e^2)^{-m_j},\quad b_j>0,\ m_j>2,    \label{eq:catalogue-energy}
\end{equation}
and take $h_j=0$ or any of
\begin{equation}
 h_j(p)=\frac{\omega_jp}{\sqrt{1+\omega_j^2p^2}},
 \qquad \tanh(\omega_jp),
 \qquad \frac2\pi\arctan(\omega_jp).            \label{eq:catalogue-odd}
\end{equation}
Then \eqref{eq:general-profile-dom} holds with $B=1$; $W,D>0$ for
$r>0$; and $Q>0$ if at least one $h_j$ is nonzero.
The logarithmic estimate \eqref{eq:general-logW} holds with
\begin{equation}
 L=2\max_{1\le j\le J}\max\{1,\nu_j\}.          \label{eq:catalogue-L}
\end{equation}
Hence Theorem \ref{thm:general-profile-class} applies on every torus
satisfying \eqref{eq:variable-torus-gap-test} for $j=1$, with any amplitude meeting
\eqref{eq:general-amplitude}.  The sum in
\eqref{eq:catalogue-mixture} is generally nonseparable in $(e,p)$ and
still describes exactly two charged species, not $2J$ species.
\end{proposition}

\begin{proof}
Every $g_j\ge1$ is even, every $h_j$ is odd with $|h_j|\le1$,
$H_j>0$, and $H_j'<0$ for $e\ge1/2$.  Positivity of the weights $a_j$
gives both domination inequalities with $B=1$.
Substitution into the moments gives
\begin{align}
 W(r)&=4\sum_{j=1}^J a_j\vartheta_j\nu_j\omega_j^2r^2
 \int_{\R^3}\frac{v_\varphi^2}{\av v}
 \frac{H_j(\av v)}
      {(1+\omega_j^2r^2v_\varphi^2)^{\nu_j+1}}\dd v>0,
                                                        \label{eq:catalogue-W}\\
 Q(r)&=2\sum_{j=1}^J a_j\int_{\R^3}
 \widehat v_\varphi[-H_j'(\av v)]h_j(rv_\varphi)\dd v,\notag\\
 D(r)&=2\sum_{j=1}^J a_j\int_{\R^3}
 [-H_j'(\av v)]g_j(rv_\varphi)\dd v>0.           \label{eq:catalogue-QD}
\end{align}
For nonzero $h_j$, the integrand defining its contribution to $Q$ is
strictly positive except on $v_\varphi=0$.
A useful upper bound, independent of the torus radii, is
\begin{equation}
 Q_\infty\le \widehat Q:=
 2\sum_{j=1}^J a_j\int_{\R^3}
 \frac{|v_\varphi|}{\av v}[-H_j'(\av v)]\dd v<\infty.
                                                        \label{eq:catalogue-Qbound}
\end{equation}
One may substitute $\widehat Q$ for $Q_\infty$ in
\eqref{eq:general-eps} to obtain a more conservative explicit threshold.

The only radial factor in each positive integrand of
\eqref{eq:catalogue-W} is
$r^2(1+\omega_j^2r^2v_\varphi^2)^{-\nu_j-1}$.
Its logarithmic derivative with respect to $\log r$ is
\[
 2-2(\nu_j+1)\frac{\omega_j^2r^2v_\varphi^2}
                        {1+\omega_j^2r^2v_\varphi^2},
\]
which belongs to $[-2\nu_j,2]$.  Differentiation under the integral,
followed by summation of positive terms, proves \eqref{eq:catalogue-L}.

For analyticity, every momentum factor extends holomorphically and
boundedly to a sufficiently narrow strip about the real axis.
The exponential energy factors have locally uniform exponential velocity
majorants after bounded potential shifts.  For the algebraic factors
the singularities in complex energy stay a fixed distance from the real
axis, and the majorants decay as $\av v^{-2m_j}$.
Cauchy estimates give the required factorial bounds for both factors.
Since $m_j>2$, even the kinetic-energy moment is finite.
These facts establish Assumption \ref{ass:analytic}, not merely
differentiability at the zero branch.
\end{proof}

The same argument allows integrals of these products against a finite
positive measure over a compact parameter set, provided the parameters
stay away from the indicated endpoints ($\beta=0$, $b=0$, $m=2$,
$\omega=0$, and so on).  Uniform strip widths and integrable majorants
then justify Fubini, differentiation, and Banach-space analyticity.
If nonzero odd factors occupy a set of positive measure, $Q>0$.
This gives an infinite-dimensional class of nonseparable profiles.

\begin{example}[Representative subclasses]
\label{ex:catalogue-subclasses}
Each line below defines a class, with the amplitude interval and the
coupling size supplied by \eqref{eq:general-endpoints} and
\eqref{eq:general-eps}.
\begin{enumerate}[label=\textup{(\alph*)},leftmargin=2.4em]
\item \emph{Arbitrary temperature and momentum scale.}
Take $J=1$, $H=e^{-\beta e}$,
$g=1+\vartheta[1-(1+\omega^2p^2)^{-\nu}]$,
$h=\tanh(\omega p)$, and $\eta(K)=K^{2n}$.
The independent parameters are $\beta,\vartheta,\omega,\nu>0$ and
$n\in\N$; the current-coupled case has $\epsilon\ne0$.
\item \emph{Algebraic-tail pair plasmas.}
Replace $H$ by $(b+e^2)^{-m}$, $b>0$, $m>2$.
The same logarithmic radial bound applies; no exponential tail is used
in the endpoint, gap, or monotonicity verification.
\item \emph{Mixed-temperature, mixed-tail populations.}
For example, with $\nu_1=\nu_2=1$, take
\begin{align*}
 M_{\rm ev}
 &=e^{-e}\left(1+\frac{p^2}{1+p^2}\right)
 +(1+e^2)^{-3}\left(1+\frac{4p^2}{1+4p^2}\right),\\
 M_{\rm od}
 &=e^{-e}\tanh p+(1+e^2)^{-3}
                 \frac{2p}{\sqrt{1+4p^2}}.
\end{align*}
Here $L=2$ and $B=1$.  Thus the very same torus
$\Omega_{1,1/5}$ from \eqref{eq:variable-torus-domain} works:
\[
 (1-a/R)^{-2L}=\frac{625}{256}
 <\frac{\ell_2}{u_1}=2.521895\ldots .
\]
This is a genuinely nonseparable common profile with nonzero coupling
for every sufficiently small $\epsilon\ne0$.
\end{enumerate}
\end{example}

\subsubsection*{Genuinely unequal, charge-balanced species}

\begin{proposition}[Unequal species with unchanged crossing coefficients]
\label{prop:balanced-species}
Assume the hypotheses of Theorem \ref{thm:general-profile-class}.  Let
$N=N(e)$ be a
nonzero analytic energy profile with the same moment bounds and
\begin{equation}
 \int_{\R^3}N(\av v)\dd v=0,\qquad
 |N(e)|\le B_NM_{\rm ev}(e,p),\qquad
 |N'(e)|\le B_N[-\partial_eM_{\rm ev}(e,p)].     \label{eq:balanced-N}
\end{equation}
The value bound holds on $\R^2$ and the derivative bound for $e\ge1/2$.
Choose $\epsilon$ satisfying the two spectral bounds in
\eqref{eq:general-eps}, and choose $\delta\ne0$ such that
\begin{equation}
 B|\epsilon|+B_N|\delta|<\tfrac12.              \label{eq:balanced-small}
\end{equation}
Then the unequal profiles
\begin{equation}
 \mu_{K,\epsilon,\delta}^{\pm}(e,p)
 =\eta(K)\bigl[M_{\rm ev}(e,p)+\epsilon M_{\rm od}(e,p)
                         \ \pm\delta N(e)\bigr]          \label{eq:balanced-profiles}
\end{equation}
are positive and energy-decreasing near the zero branch, and have
$A=\phi=0$ as an exact equilibrium for every $K$.
Their coefficients $c_{ij}$ along that branch are still
\eqref{eq:general-coefficients}.  In particular the same $K_*$,
transversality, and \textup{(C1)--(C4)} apply, with genuine two-field
coupling if $\epsilon\ne0$ and $Q>0$.
\end{proposition}

\begin{proof}
The smallness bound proves positivity and energy monotonicity.
At the zero branch the charge difference is $2\eta\delta N(\av v)$;
its integral is zero by \eqref{eq:balanced-N}, and its current moment
vanishes by oddness of $\widehat v_\varphi$.
Every $c_{ij}$ is a sum of the two species' derivatives, so the
$+\delta N$ and $-\delta N$ contributions cancel.  This proves every
claim about the stationary linearization.
\end{proof}

There are explicit choices of $N$.  In the two-term mixture in
Example \ref{ex:catalogue-subclasses}\textup{(c)}, put
\[
 H_1(e)=e^{-e},\quad H_2(e)=(1+e^2)^{-3},\quad
 n_j=\int_{\R^3}H_j(\av v)\dd v,\quad
 c=n_1/n_2,\quad N=H_1-cH_2.
\]
Then $n_j$ are finite positive one-dimensional radial integrals,
$\int N(\av v)\dd v=0$, and \eqref{eq:balanced-N} holds with
$B_N=\max\{1,c\}$ since $g_j\ge1$ and $H_j'<0$ for positive energies.
Thus $(|\epsilon|+\max\{1,c\}|\delta|)<1/2$, together with the spectral
coupling bounds, gives a fully specified unequal-species class on
$\Omega_{1,1/5}$.  This extension changes the nonlinear branch coefficients
in general, although it leaves its linear bifurcation point unchanged.

\subsubsection*{Location of the crossing and limitations}

For the equal or charge-balanced classes, write $t=\eta(K)$ and
\[
 S_\epsilon(t)=L_A-tW+
       \epsilon^2t^2Q(L_\phi+tD)^{-1}Q.
\]
If $t_*(\epsilon)=\eta(K_*(\epsilon))$ is the crossing supplied above,
the min--max principle immediately gives
\begin{equation}
 \Lambda_1(W)\le t_*(\epsilon)\le t_+,\qquad
 \frac{\lambda_{A,1}}{W_+}\le\Lambda_1(W)
                    \le\frac{\lambda_{A,1}}{W_-}.       \label{eq:general-location}
\end{equation}
The first inequality is strict when $\epsilon\ne0$ and $Q>0$.
Indeed, $L_A-\Lambda_1(W)W$ is nonnegative with simple kernel and the
Schur correction has positive quadratic form on that kernel.
Compact resolvent excludes a zero eigenvalue after this correction.

More precisely, let $\psi_0>0$ solve
$L_A\psi_0=t_0W\psi_0$, with $t_0=\Lambda_1(W)$ and
$\int_\Omega W\psi_0^2\dd x=1$.  Analytic perturbation of its simple
zero eigenvalue \cite{Kato} gives
\begin{equation}
 t_*(\epsilon)=t_0+\epsilon^2t_0^2
 \left\langle Q\psi_0,(L_\phi+t_0D)^{-1}Q\psi_0\right\rangle
 +O(\epsilon^4).                               \label{eq:general-shift}
\end{equation}
To verify the coefficient, set $\xi=\epsilon^2$, differentiate
$S_{\sqrt{\xi}}(t(\xi))a(\xi)=0$ at $\xi=0$, and pair with $\psi_0$.
The term containing $a'(0)$ vanishes and the normalization makes the
coefficient of $-t'(0)$ equal to one.  The shift in $K_*$ is the
coefficient in \eqref{eq:general-shift} divided by
$\eta'(K_*(0))$.  This is a shift of the \emph{stationary} singularity;
it is not a statement about dynamical stability.

The strict inequalities above also persist under sufficiently small
changes of the moments $W,Q,D$ and amplitude in the indicated norms,
provided the exact reference-branch cancellation is retained.
Arbitrary independent changes in the two species need not retain that
branch, and are not covered by this persistence statement.
Likewise, neither an arbitrary thick torus nor an amplitude of
insufficient range is covered.  The spectral ratio
\eqref{eq:general-window} is sufficient, not necessary.

The generalization uses the classical min--max and simple-eigenvalue
bifurcation tools \cite{Kato,Amann,CrandallRabinowitz1971,Rabinowitz1971}.
It is a profile-level verification theorem for the present stationary
RVM reduction, not a new abstract bifurcation theorem, as distinguished
in the introduction.

\section{Examples of norm-unbounded two-species RVM curves}
\label{sec:global-bif}

This section gives electric-only, magnetic-only, and simultaneous
norm-escape constructions of RVM equilibria solution curves.  The electric-only class is an unequal-activity
Maxwell--J\"uttner pair for which \(A\equiv0\) and only the electric
potential escapes.  Equal-species even-momentum profiles with
\(\phi\equiv0\) give three distinct magnetic-only mechanisms: a
diverging current coefficient, finite-parameter asymptotic resonance,
and a superlinear positive-cone continuum at bounded parameter.  The
simultaneous class combines activity imbalance with current generation
and forces both potentials to diverge along one positive global continuum.
The calculations identify every coefficient from velocity moments.
Assume throughout that $\Omega$ is a connected smooth solid torus with
$0<r_-\le r\le r_+$.

For positive-cone continuation in this section, the trace-zero space and
its nonnegative cone are
\[
Y=\{u\in C^1_{\rm ax}(\overline\Omega):u|_{\partial\Omega}=0\},
\qquad Y_+=\{u\in Y:u\ge0\}.
\]
We do not mean the $C^1$-closure of compactly supported functions.
The cone interior consists of functions $u>0$ in $\Omega$ with
$\partial_nu<0$ on $\partial\Omega$.  For fixed $p>3$, bounded
right-hand sides have bounded Dirichlet $W^{2,p}$ norms and are
therefore precompact in $Y$; the same estimate applies to normalized
nonreturn sequences in the proofs below.

\subsection{Case E: unbounded electric potential and bounded magnetic
potential}

Put
\begin{equation}
 m_0=\int_{\R^3}e^{-\av v}\dd v,
 \qquad L_\phi\varphi_1=\lambda_{\phi,1}\varphi_1,
 \qquad \varphi_1>0.
 \label{eq:electric-juttner-data}
\end{equation}
Let \(a_\pm:I\to(0,\infty)\) be real analytic and take the
two-species profiles
\begin{equation}
 \mu_K^+(e,p)=a_+(K)e^{-e},\qquad
 \mu_K^-(e,p)=a_-(K)e^{-e}.                       \label{eq:electric-juttner-profiles}
\end{equation}
They are Maxwell--J\"uttner distributions at the normalized inverse
temperature one, with activities \(a_\pm\).  Since they are independent
of canonical momentum, their azimuthal current is zero and hence
\(A_K\equiv0\).  The remaining equation is the relativistic
Poisson--Boltzmann problem
\begin{equation}
 L_\phi\phi=m_0\bigl(a_+(K)e^{-\phi}
                         -a_-(K)e^\phi\bigr),\qquad
 \phi|_{\partial\Omega}=0.                        \label{eq:electric-juttner-PB}
\end{equation}

\begin{theorem}[Electric-only J\"uttner norm escape]
\label{thm:electric-juttner-escape}
For every \(K\in I\), equation \eqref{eq:electric-juttner-PB} has a
unique solution \(\phi_K\in X_\phi\), and \(K\mapsto\phi_K\) is real
analytic.  If \(a_+(K)\ge a_-(K)\) and
\begin{equation}
 \delta(K)=\frac12\log\frac{a_+(K)}{a_-(K)},       \label{eq:electric-delta}
\end{equation}
then \(0\le\phi_K\le\delta(K)\).  For each prescribed \(R>0\),
\begin{equation}
 \boxed{\quad
 m_0\bigl(a_+(K)e^{-R}-a_-(K)e^R\bigr)
   >\lambda_{\phi,1}R
 \quad\Longrightarrow\quad \|\phi_K\|_\infty>R .\quad}       \label{eq:electric-juttner-test}
\end{equation}
In particular, for
\begin{equation}
 a_+(K)=e^{\gamma K},\qquad a_-(K)=e^{-\gamma K},
 \qquad \gamma>0,\qquad K\ge0,                 \label{eq:electric-activities}
\end{equation}
the sufficient condition becomes
\begin{equation}
 2m_0\sinh(\gamma K-R)>\lambda_{\phi,1}R,          \label{eq:electric-explicit-test}
\end{equation}
and
\begin{equation}
 \|\phi_K\|_{L^\infty}+\|\phi_K\|_{X_\phi}\longrightarrow\infty,
 \qquad A_K\equiv0,\qquad K\to\infty.          \label{eq:electric-only-divergence}
\end{equation}
Thus \(K\mapsto(0,\phi_K,K)\) is a global analytic continuation curve
with electric-only norm escape.
\end{theorem}

\begin{proof}
Fix $K$ and write $F(s)=m_0(a_+e^{-s}-a_-e^s)$ and
$\delta=\frac12\log(a_+/a_-)$.  The constants
$m=\min\{0,\delta\}$ and $M=\max\{0,\delta\}$ are respectively a
subsolution and a supersolution, since $F(m)\ge0$ and $F(M)\le0$.
Truncate $F$ outside $[m,M]$ to the bounded continuous function
$F_T(s)=F(\min\{M,\max\{m,s\}\})$.  The compact Dirichlet
resolvent and Schauder's fixed-point theorem give a solution of
$L_\phi\phi=F_T(\phi)$.  Testing above $M$ and below $m$ with the
positive and negative parts, respectively, shows $m\le\phi\le M$;
hence it solves the untruncated equation.  Since $F'(s)<0$, the
difference of any two solutions tested against itself vanishes only
when the solutions coincide.  Elliptic regularity gives $\phi\in
X_\phi$.  This barrier argument avoids presuming that the exponential
functional is finite on every element of three-dimensional $H^1_0$.
The field derivative is
\[
 L_\phi+m_0\bigl(a_+(K)e^{-\phi_K}
                  +a_-(K)e^{\phi_K}\bigr)>0,
\]
so the analytic implicit-function theorem gives the global analytic
graph.  The ordered-barrier argument above gives
$\min\{0,\delta(K)\}\le\phi_K\le\max\{0,\delta(K)\}$.
When $a_+(K)\ge a_-(K)$, one has $\delta(K)\ge0$ and hence
$0\le\phi_K\le\delta(K)$.

If \(\|\phi_K\|_\infty\le R\), pairing
\eqref{eq:electric-juttner-PB} with \(\varphi_1\) yields
\[
 \lambda_{\phi,1}\int_\Omega\varphi_1\phi_K\dd x
 \ge m_0\bigl(a_+(K)e^{-R}-a_-(K)e^R\bigr)
        \int_\Omega\varphi_1\dd x.
\]
The left side is at most
\(\lambda_{\phi,1}R\int\varphi_1\), proving
\eqref{eq:electric-juttner-test}.  Substitution of
\eqref{eq:electric-activities} gives
\eqref{eq:electric-explicit-test}, which holds for every fixed \(R\)
when \(K\) is large.  The embedding
\(X_\phi\hookrightarrow C^0(\overline\Omega)\) proves
\eqref{eq:electric-only-divergence}.
\end{proof}

This example is a regular continuation rather than an electric
secondary bifurcation.  That distinction is forced by physics: for
energy-decreasing profiles the electric block is strictly positive, so
loss of invertibility must enter through the magnetic Schur complement.
The activities in \eqref{eq:electric-juttner-profiles} are the
relativistic analogues of Boltzmann fugacities; the resulting
charge-separation equation is the Poisson--Boltzmann response of a
non-neutral thermal pair plasma \cite{Juttner1911,Davidson2001}.

\subsection{Case M: bounded electric potential and unbounded magnetic
potential}

The first magnetic class has a diverging current coefficient and a
unique positive solution for every parameter beyond its crossing.
The second escapes at a finite asymptotic resonance.  Both have
decreasing response quotient \(H(r,s)/s\).  The third has an
increasing response quotient; Proposition
\ref{prop:superlinear-case-M} places the counter-rotating
Maxwell--J\"uttner pitchfork of Example
\ref{ex:counterrotating-juttner-pitchfork}.  In all three
classes \(\phi\equiv0\).

\subsubsection{Sublinear coefficient growth.}

Let $g:\R\to(0,\infty)$ be even, real analytic,
integrable, strictly decreasing on $(0,\infty)$, and suppose its
derivatives have velocity-integrable analytic majorants.  Put
\begin{equation}
 G(p)=\int_0^p(p-t)g(t)\dd t.                       \label{eq:unbounded-G}
\end{equation}
Then $G$ is even and nonnegative, $G(0)=G'(0)=0$, and $G''=g$.
Let $C>0$, let $\beta:\R\to(0,\infty)$ be real analytic, and define
\begin{equation}
 \mu_K^+(e,p)=\mu_K^-(e,p)
 =e^{-e}\bigl(C+\beta(K)G(p)\bigr).                \label{eq:unbounded-profiles}
\end{equation}
Both species are strictly positive and energy decreasing.
This is the neutral magnetic pitchfork mechanism:
the $p$-independent background $Ce^{-e}$ carries no current, while
$M(e,p)=e^{-e}G(p)$ and $\eta(K)=\beta(K)$ generate the crossing.
In particular $(A,\phi)=(0,0)$ is an exact reference branch.

For $r\in[r_-,r_+]$ set
\begin{align}
 h_r(t)&=\int_{\R^3}e^{-\av v}g\bigl(r(v_\varphi+t)\bigr)\dd v,
                                                        \label{eq:unbounded-h}\\
 H(r,s)&=r^2\int_{-s}^{s}h_r(t)\dd t,\qquad
 \alpha(r)=2r^2h_r(0),                               \label{eq:unbounded-H}\\
 \Lambda&=\Lambda_1(\alpha)
 =\inf_{a\ne0}
 \frac{\displaystyle\int_\Omega(|\nabla a|^2+r^{-2}a^2)\dd x}
      {\displaystyle\int_\Omega\alpha a^2\dd x}.  \label{eq:unbounded-Lambda}
\end{align}
Parity makes $\phi=0$ invariant.  Integration by parts in
$v_\varphi$ gives the exact equation
\begin{equation}
 L_AA=\beta(K)H(r,A),\qquad A|_{\partial\Omega}=0.  \label{eq:unbounded-scalar}
\end{equation}
The convolution/layer-cake argument used below shows that $h_r$ is even
and strictly decreasing on $(0,\infty)$.  Hence
\begin{equation}
 s\longmapsto\frac{H(r,s)}s
 =2r^2\frac1s\int_0^sh_r(t)\dd t                   \label{eq:unbounded-sublinear}
\end{equation}
is strictly decreasing from $\alpha(r)$ to zero, uniformly in
$r\in[r_-,r_+]$.

For $R>0$ define the computable finite-amplitude quotient
\begin{equation}
 q_R=\min_{\substack{r_-\le r\le r_+\\0\le s\le R}}
 Q(r,s),\qquad
 Q(r,s)=\frac{H(r,s)}{\alpha(r)s}\ (s>0),\qquad
 Q(r,0)=1.                                         \label{eq:unbounded-qR}
\end{equation}
The value at $s=0$ is the continuous extension.  Hence
$0<q_R\le1$.

\begin{theorem}[A computable unbounded RVM bifurcation curve]
\label{thm:explicit-unbounded-rvm}
Assume
\begin{equation}
 \beta'(K)>0,\qquad
 \lim_{K\to-\infty}\beta(K)<\Lambda,\qquad
 \lim_{K\to+\infty}\beta(K)=\infty.               \label{eq:unbounded-beta}
\end{equation}
There is a unique $K_*$ with $\beta(K_*)=\Lambda$.  The zero equilibrium
at $K_*$ is a simple transverse bifurcation point.  For every $K>K_*$,
\eqref{eq:unbounded-scalar} has a unique positive solution $A_K$ and a
unique negative solution $-A_K$; there is no nonzero one-signed solution
for $K\le K_*$.  The maps $K\mapsto\pm A_K$ are real analytic,
\begin{equation}
 \|A_K\|_{X_A}\longrightarrow0\quad(K\downarrow K_*),
 \qquad
 \|A_K\|_{L^\infty}+\|A_K\|_{X_A}\longrightarrow\infty
 \quad(K\to\infty),                                \label{eq:unbounded-arm-limits}
\end{equation}
and the two arms lie in the Rabinowitz component issuing from
$(0,0,K_*)$.  In particular that component is unbounded.
The local pitchfork has the explicit direction
\begin{equation}
 A(s)=s\psi+O(s^3),\qquad
 K(s)=K_*+\kappa_2s^2+O(s^4),\qquad
 \kappa_2=-\frac{\Lambda}{3\beta'(K_*)}
 \frac{\displaystyle\int_\Omega r^2h_r''(0)\psi^4\dd x}
      {\displaystyle\int_\Omega\alpha\psi^2\dd x}>0,
 \label{eq:unbounded-local-pitchfork}
\end{equation}
where $\psi>0$ solves $L_A\psi=\Lambda\alpha\psi$.

More quantitatively, for every $R>0$,
\begin{equation}
 \boxed{\quad
 \beta(K)>\frac{\Lambda}{q_R}
 \quad\Longrightarrow\quad
 \|A_K\|_{L^\infty(\Omega)}>R .\quad}              \label{eq:unbounded-growth-test}
\end{equation}
\end{theorem}

\begin{proof}
The linearization on the magnetic invariant subspace is
$L_A-\beta(K)\alpha$.  Its principal kernel at $K_*$ is generated by a
positive eigenfunction $\psi$, and
\[
 \left\langle\psi,
 \partial_K(L_A-\beta(K)\alpha)|_{K=K_*}\psi\right\rangle
 =-\beta'(K_*)\int_\Omega\alpha\psi^2\dd x<0.
\]
The electric block is positive by energy monotonicity, and parity makes
the two off-diagonal blocks vanish.  Thus this is a simple transverse
crossing for the full two-field system.
Moreover,
\[
 h_r''(0)=r\int_{\R^3}
 e^{-\av v}\widehat v_\varphi g'(rv_\varphi)\dd v<0
 \quad(r>0).
\]
Since $H(r,A)=\alpha(r)A+\frac13r^2h_r''(0)A^3+O(A^5)$,
projection of the cubic equation onto $\psi$ gives
\eqref{eq:unbounded-local-pitchfork}.  This also proves that the
global arms are the continuation of the local neutral pitchfork.

The sublinear positive-solution theorem of Brezis and Oswald
\cite{BrezisOswald1986}, applied to \eqref{eq:unbounded-scalar}, gives a
positive solution precisely when $\beta(K)>\Lambda$, and gives
uniqueness.  Oddness supplies the negative solution.  Moreover,
\eqref{eq:unbounded-sublinear} implies
\[
 \bigl(L_A-\beta(K)\partial_sH(r,A_K)\bigr)A_K
 =\beta(K)\bigl[H(r,A_K)-A_K\partial_sH(r,A_K)\bigr]>0.
\]
The magnetic linearization therefore has positive principal eigenvalue.
Along the entire invariant set \(\phi=0\), equality of the two species
and evenness in \(p\) give \(c_{12}=c_{21}=0\), while energy
monotonicity keeps the electric block positive.  Hence the full
derivative is invertible along both arms.  The analytic
implicit-function theorem, local Crandall--Rabinowitz description, and
compactness give the asserted analytic continuation and the first limit
in \eqref{eq:unbounded-arm-limits}.

If $\|A_K\|_\infty\le R$, definition \eqref{eq:unbounded-qR} and pairing
with $\psi$ give
\[
 \Lambda\int_\Omega\alpha\psi A_K\dd x
 =\beta(K)\int_\Omega\psi H(r,A_K)\dd x
 \ge\beta(K)q_R\int_\Omega\alpha\psi A_K\dd x.
\]
Thus $\beta(K)\le\Lambda/q_R$, proving
\eqref{eq:unbounded-growth-test}.  Since $R$ is arbitrary and
$\beta(K)\to\infty$, the $L^\infty$ norm diverges.  The embedding
$X_A\hookrightarrow C^0(\overline\Omega)$ gives the $X_A$ divergence.
Each arm meets the reference branch at the Crandall--Rabinowitz point,
so it belongs to the corresponding global component.
\end{proof}

The implication \eqref{eq:unbounded-growth-test} certifies a point of
arbitrarily large field norm directly from \(\beta(K)\) and \(q_R\).

\begin{example}[The logarithmic-cosh family]
\label{ex:logcosh-unbounded}
For \(g(p)=\operatorname{sech}^2p\), \(G(p)=\log\cosh p\), and
\(\beta(K)=\Lambda e^{\sigma(K-K_*)}\) with \(\sigma>0\), the
pitchfork of Theorem \ref{thm:explicit-unbounded-rvm} is supercritical.
The shift bound \(\cosh(x+y)\le e^{|y|}\cosh x\) gives
\begin{equation}
 q_R\ge\frac{1-e^{-2r_+R}}{2r_+R},\qquad
 \beta(K)>\frac{2r_+R\Lambda}{1-e^{-2r_+R}}
 \ \Longrightarrow\ \|A_K\|_\infty>R.
 \label{eq:unbounded-qR-explicit}
\end{equation}
A physical-ray choice is \(\beta(K)=K\), \(K\ge0\), with
\(K_*=\Lambda\) and \(\mu_K^\pm=e^{-e}(C+K\log\cosh p)\).
On \(\Omega_{1,1/5}\), write
\(\alpha_-=\min_{[4/5,6/5]}\alpha\) and
\(\alpha_+=\max_{[4/5,6/5]}\alpha\), and use the eigenvalue
enclosures \(\ell_1,u_1\) of \eqref{eq:variable-torus-mode-enclosures}:
\begin{equation}
 \frac{\ell_1}{\alpha_+}\le\Lambda\le
 \frac{u_1}{\alpha_-},\qquad
 K>\frac{u_1}{\alpha_-}
 \frac{(12/5)R}{1-e^{-(12/5)R}}
 \ \Longrightarrow\ \|A_K\|_\infty>R.
 \label{eq:unbounded-torus-test}
\end{equation}
Only the one-dimensional velocity moment in \(\alpha_\pm\)
requires numerical enclosure.
\end{example}

\begin{example}[A general algebraic-curvature class]
\label{ex:algebraic-curvature-unbounded}
For $q>1/2$, let
\[
 g_q(p)=(1+p^2)^{-q},\qquad
 G_q(p)=\int_0^p(p-t)(1+t^2)^{-q}\dd t.
\]
Then \eqref{eq:unbounded-profiles} satisfies every hypothesis of Theorem
\ref{thm:explicit-unbounded-rvm}.  Since
$|(\log g_q)'|\le q$, one has the explicit certificate
\begin{equation}
 q_R\ge\frac{1-e^{-qr_+R}}{qr_+R}.                 \label{eq:algebraic-qR}
\end{equation}
Two elementary members are
\[
 G_1(p)=p\arctan p-\tfrac12\log(1+p^2),\qquad
 G_{3/2}(p)=\sqrt{1+p^2}-1.
\]
Thus the construction is an infinite-dimensional profile class, not a
single tuned example.
\end{example}

For the logarithmic-cosh family, the source also satisfies
\begin{equation}
 |H(r,s)|\le 2r\int_{\R^3}e^{-\av v}\dd v.         \label{eq:unbounded-source-upper}
\end{equation}
Hence $\|A_K\|_{X_A}\le C\beta(K)$ on every finite parameter interval,
the profiles never lose admissibility, and no finite-$K$ blowup occurs.
The branch therefore escapes simultaneously through
$K\to\infty$ and $\|A_K\|_{X_A}\to\infty$.

\subsubsection{Finite-parameter asymptotic resonance.}

The integrability of \(g\) above makes \(H(r,s)/s\to0\); consequently
large fields require \(\beta(K)\to\infty\).  A second, equally
computable mechanism appears when the momentum curvature has a positive
limit.  Fix \(g_\infty>0\), let \(\widetilde g:\R\to(0,\infty)\) satisfy
the analyticity, evenness, integrability, and strict-decrease assumptions
imposed on \(g\) above, and set
\begin{equation}
 g(p)=g_\infty+\widetilde g(p),\qquad
 G(p)=\frac{g_\infty}{2}p^2+\widetilde G(p),\qquad
 \widetilde G(p)=\int_0^p(p-t)\widetilde g(t)\dd t.
 \label{eq:resonant-curvature}
\end{equation}
For \(K>0\), define the equal two-species profiles
\begin{equation}
 \mu_K^+(e,p)=\mu_K^-(e,p)
 =e^{-e}\bigl(C+K G(p)\bigr),\qquad C>0.           \label{eq:resonant-profiles}
\end{equation}
They are strictly positive, energy decreasing, and satisfy Assumption
\ref{ass:analytic}; the quadratic growth of \(G\) is harmless because of
the exponential energy factor.

Put
\begin{align}
 m_0&=\int_{\R^3}e^{-\av v}\dd v,
 &\widetilde h_r(t)&=\int_{\R^3}e^{-\av v}
     \widetilde g\bigl(r(v_\varphi+t)\bigr)\dd v,\notag\\
 \alpha_\infty(r)&=2r^2g_\infty m_0,
 &\widetilde H(r,s)&=r^2\int_{-s}^{s}\widetilde h_r(t)\dd t,\notag\\
 H_\infty(r,s)&=\alpha_\infty(r)s+\widetilde H(r,s),
 &\alpha_0(r)&=\alpha_\infty(r)+2r^2\widetilde h_r(0).
 \label{eq:resonant-coefficients}
\end{align}
Then the invariant equation on \(\phi=0\) is
\begin{equation}
 L_AA=K H_\infty(r,A).                             \label{eq:resonant-scalar}
\end{equation}
The quotient \(H_\infty(r,s)/s\) is strictly decreasing from
\(\alpha_0(r)\) to \(\alpha_\infty(r)\), uniformly in
\(r\in[r_-,r_+]\).  Define the two computable weighted principal
eigenvalues
\begin{equation}
 \Lambda_0=\Lambda_1(\alpha_0),\qquad
 \Lambda_\infty=\Lambda_1(\alpha_\infty).
 \label{eq:resonant-eigenvalues}
\end{equation}
Since \(\alpha_0>\alpha_\infty\), strict monotonicity of the principal
Rayleigh quotient gives \(0<\Lambda_0<\Lambda_\infty\).

For \(R>0\), let
\begin{equation}
 d_R=\min_{\substack{r_-\le r\le r_+\\0\le s\le R}}D(r,s),
 \quad
 D(r,s)=\frac{\widetilde H(r,s)}{\alpha_\infty(r)s}\ (s>0),
 \quad
 D(r,0)=\frac{2r^2\widetilde h_r(0)}{\alpha_\infty(r)}.
 \label{eq:resonant-dR}
\end{equation}
This is a minimum of a positive continuous function on a compact
rectangle, so \(d_R>0\), and it is rigorously computable by velocity
tail estimates and interval quadrature.

\begin{theorem}[Finite-parameter norm-unbounded RVM curve]
\label{thm:finite-resonance-rvm}
For the profiles \eqref{eq:resonant-profiles}, the neutral equilibrium
at \(K=\Lambda_0\) is a simple transverse bifurcation point.  For every
\begin{equation}
 \Lambda_0<K<\Lambda_\infty                       \label{eq:resonant-window}
\end{equation}
there is a unique positive solution \(A_K\) of
\eqref{eq:resonant-scalar} and a unique negative solution \(-A_K\).
There is no nonzero one-signed solution for
\(K\le\Lambda_0\) or \(K\ge\Lambda_\infty\).  Both arms are real
analytic, lie in the global component issuing from
\((0,0,\Lambda_0)\), and satisfy
\begin{equation}
 \|A_K\|_{X_A}\longrightarrow0\quad(K\downarrow\Lambda_0),
 \qquad
 \|A_K\|_{L^\infty}+\|A_K\|_{X_A}\longrightarrow\infty
 \quad(K\uparrow\Lambda_\infty).                 \label{eq:resonant-limits}
\end{equation}
For every prescribed \(R>0\), the computable implication
\begin{equation}
 \boxed{\quad
 K>\frac{\Lambda_\infty}{1+d_R}
 \quad\Longrightarrow\quad \|A_K\|_{L^\infty(\Omega)}>R
 \quad}                                           \label{eq:resonance-growth-test}
\end{equation}
holds.  Thus \(K=\Lambda_\infty\) is an asymptotic bifurcation value,
or bifurcation point from infinity, while the physical parameter remains
bounded.
\end{theorem}

\begin{proof}
The linearization at zero is \(L_A-K\alpha_0\); hence the principal
kernel at \(K=\Lambda_0\) is simple, and its crossing form is
\(-\int_\Omega\alpha_0\psi_0^2\dd x<0\).  The electric block is
positive and the off-diagonal blocks vanish by parity, so the crossing
is simple also in the full two-field system.

For \(s>0\), the quotient \(H_\infty(r,s)/s\) is strictly decreasing,
with limits \(\alpha_0\) at zero and \(\alpha_\infty\) at infinity.
The positive-solution criterion of Brezis--Oswald
\cite{BrezisOswald1986}, equivalently its compact positivity-improving
resolvent form \cite{Klimsiak2026}, therefore gives existence and
uniqueness precisely in the spectral window
\eqref{eq:resonant-window}.  For completeness, necessity follows
directly.  If \(\psi_0>0\) and \(\psi_\infty>0\) are the principal
eigenfunctions for \(\alpha_0\) and \(\alpha_\infty\), respectively,
then any positive solution obeys
\begin{align*}
 \Lambda_0\int_\Omega\alpha_0\psi_0A\dd x
  &=K\int_\Omega\psi_0H_\infty(r,A)\dd x
    <K\int_\Omega\alpha_0\psi_0A\dd x,\\
 \Lambda_\infty\int_\Omega\alpha_\infty\psi_\infty A\dd x
  &=K\int_\Omega\psi_\infty H_\infty(r,A)\dd x
    >K\int_\Omega\alpha_\infty\psi_\infty A\dd x.
\end{align*}
Thus \(K>\Lambda_0\) and \(K<\Lambda_\infty\).

Strict decrease of the quotient also gives
\[
 \bigl(L_A-K\partial_sH_\infty(r,A_K)\bigr)A_K
 =K\bigl[H_\infty(r,A_K)
          -A_K\partial_sH_\infty(r,A_K)\bigr]>0.
\]
The magnetic derivative has positive principal eigenvalue.  Equality of
the species and evenness of \(G\) make \(c_{12}=c_{21}=0\) at every
point with \(\phi=0\), and the electric block is positive by energy
monotonicity.  Thus the full derivative is invertible along both
nontrivial arms.  Analyticity and the
local simple-eigenvalue theorem therefore give analytic global graphs.
Uniqueness identifies their lower limits with the local branch, so
\(A_K\to0\) as \(K\downarrow\Lambda_0\).

If \(\|A_K\|_\infty\le R\), testing with \(\psi_\infty\) and using
\eqref{eq:resonant-dR} yields
\begin{align*}
 \Lambda_\infty\int_\Omega\alpha_\infty\psi_\infty A_K\dd x
 &=K\int_\Omega\psi_\infty
       \bigl(\alpha_\infty A_K+\widetilde H(r,A_K)\bigr)\dd x\\
 &\ge K(1+d_R)
       \int_\Omega\alpha_\infty\psi_\infty A_K\dd x.
\end{align*}
This proves \eqref{eq:resonance-growth-test}.  Because its threshold is
strictly below \(\Lambda_\infty\) for every fixed \(R\), the
\(L^\infty\) norm tends to infinity as
\(K\uparrow\Lambda_\infty\); the embedding of \(X_A\) into
\(C^0(\overline\Omega)\) gives the \(X_A\)-norm divergence.  Oddness
of \(H_\infty\) supplies the negative arm.
\end{proof}

\begin{example}[Quadratic plus logarithmic-cosh profiles]
\label{ex:finite-resonance-logcosh}
For \(\widetilde g=\operatorname{sech}^2p\), take
\begin{equation}
 \mu_K^\pm(e,p)=e^{-e}
 \left(C+K\left[\tfrac12g_\infty p^2+\log\cosh p\right]\right),
 \qquad g_\infty>0.                       \label{eq:finite-resonance-logcosh-profile}
\end{equation}
With
\(\eta_-=\min_{[r_-,r_+]}\widetilde h_r(0)/(g_\infty m_0)>0\),
the same shift bound as \eqref{eq:unbounded-qR-explicit} yields
\begin{equation}
 d_R\ge\underline d_R:=
 \eta_-\frac{1-e^{-2r_+R}}{2r_+R},\qquad
 \frac{\Lambda_\infty}{1+\underline d_R}<K<\Lambda_\infty
 \ \Longrightarrow\ \|A_K\|_\infty>R.
 \label{eq:resonant-dR-explicit}
\end{equation}
The interval is nonempty for every \(R\).  On a circular torus,
Proposition \ref{prop:general-torus} bounds \(\Lambda_\infty\)
by Bessel-zero and moment enclosures, making this a finite test.
\end{example}


\subsubsection{Superlinear thermal response.}

\begin{proposition}[Superlinear magnetic norm escape]
\label{prop:superlinear-case-M}
Suppose the invariant neutral sector of a two-species family has
\[
 L_AA=K F(r,A),\qquad K\ge0,\qquad \phi=0,
\]
where \(F\) is odd and locally smooth, \(F(r,s)>0\) for \(s>0\), and,
uniformly for \(r\in[r_-,r_+]\),
\begin{equation}
 \partial_sF(r,0)=w(r)>0,\qquad F(r,s)>w(r)s
 \quad(s>0),\qquad K_*:=\Lambda_1(w)>0.
 \label{eq:superlinear-case-M}
\end{equation}
Assume the Dirichlet resolvent is compact and strongly positive,
and the electric block is invertible at \((K_*,0,0)\).  The local
positive branch belongs to a connected continuum \(\mathscr C_+^M\)
such that
\begin{equation}
 0<K<K_*\quad\hbox{at every nonzero point},\qquad
 \sup_{(K,A)\in\mathscr C_+^M}\|A\|_\infty=\infty.
 \label{eq:superlinear-case-M-escape}
\end{equation}
For any \(R>0\) the finite-amplitude number
\begin{equation}
 M_R=\max_{\substack{r_-\le r\le r_+\\0\le s\le R}}
 \frac{F(r,s)}{w(r)s},\qquad
 \left.\frac{F(r,s)}{w(r)s}\right|_{s=0}:=1,
 \label{eq:superlinear-MR}
\end{equation}
gives the computable test
\begin{equation}
 0<K<\frac{K_*}{M_R}\quad\Longrightarrow\quad
 \|A\|_\infty>R
 \quad\hbox{for every nonzero positive solution}.
 \label{eq:superlinear-case-M-test}
\end{equation}
\end{proposition}

\begin{proof}
Let \(\psi>0\) solve \(L_A\psi=K_*w\psi\).  Testing the scalar
equation against \(\psi\) gives
\[
 (K_*-K)\int_\Omega w\psi A\dd x
 =K\int_\Omega\psi\,[F(r,A)-wA]\dd x>0.
\]
Work in the trace-zero space $Y$ and cone $Y_+$ defined at the start
of this section.  The compact Dirichlet map
$A\mapsto K L_A^{-1}F(r,A)$ is strongly positive for
$K>0$ and $A\in Y_+\setminus\{0\}$; at $K=0$ it is the zero map.
The strong maximum principle and Hopf lemma
exclude contact with $\partial Y_+$ away from zero.  At $K=0$ the
only solution is zero.  If positive solutions tended to a zero
state at $K_1$, normalize them in $Y$ and use elliptic compactness:
the nonzero nonnegative limit solves $L_Aa=K_1wa$ and is the
positive principal eigenfunction, so $K_1=K_*$.  The strict test above
excludes nonzero positive solutions at $K=K_*$; the local positive
germ is the unique one at that point.  Thus the unilateral
positive-cone alternative \cite{DancerPositive1973} has no return or
cone-boundary alternative in the strip $0\le K\le K_*$.
If $K$ and $\|A\|_\infty$ were bounded, local boundedness of $F$
and elliptic estimates would bound $\|A\|_{X_A}$ and make the zero
set compact.  The continuum must therefore escape in field norm,
proving \eqref{eq:superlinear-case-M-escape}.
If \(0<A\le R\), then \(F(r,A)\le M_Rw(r)A\); the same eigenfunction
test gives \(K_*\le KM_R\), proving
\eqref{eq:superlinear-case-M-test}.
\end{proof}

This criterion requires neither decreasing \(F(r,s)/s\), uniqueness
of positive solutions, nor a global graph.  It does not imply
\(K\downarrow0\) at the unbounded end of \(\mathscr C_+^M\): blow-up at
another parameter in \((0,K_*)\) remains possible.

The next example uses a standard kinetic thermal profile, rather than
the curvature-designed functions above.  It connects a zero-current
counter-rotating family to two current-carrying families by a
subcritical pitchfork and satisfies
\eqref{eq:superlinear-case-M} exactly.

\begin{example}[Counter-rotating thermal pitchfork]
\label{ex:counterrotating-juttner-pitchfork}
On a smooth solid torus fix \(C,\vartheta,\omega>0\) with
\(\omega r_+<\vartheta\), and set
\begin{equation}
 \mu_K^+(e,p)=\mu_K^-(e,p)
 =e^{-\vartheta e}\bigl(C+K\cosh(\omega p)\bigr),\qquad K\ge0.
 \label{eq:rotor-juttner-profile}
\end{equation}
This is a positive, energy-decreasing equal mixture of oppositely
rotating thermal populations; the light-cylinder bound gives
integrable analytic moments \cite{Juttner1911,deGrootEtAl1980,Davidson2001}.
Canonical-momentum exponentials also occur in Harris--Channell and
rigid-rotor equilibrium models
\cite{Harris1962,Channell1976,ChenDavidsonPakter1997,GedalinEtAl2011}.
With \(q=\omega r\) and
\(\Delta(q)=\sqrt{\vartheta^2-q^2}\), the needed moments are
\begin{equation}
 Z_\vartheta(q)=\frac{4\pi\vartheta}{\Delta(q)^2}K_2(\Delta(q)),
 \quad m_\vartheta=Z_\vartheta(0),\quad
 I_\omega(r)=\frac{4\pi q}{\Delta(q)^2}K_2(\Delta(q))>0,
 \label{eq:rotor-I-bessel}
\end{equation}
where \(K_2\) is the modified Bessel function; these identities
follow by Lorentz-invariant integration
\cite{NISTDLMF,Juttner1911}.  Direct velocity integration gives the
full charge and current before neutrality is imposed:
\begin{equation}
 \rho_K=-2\bigl[Cm_\vartheta+
 KZ_\vartheta(q)\cosh(qA)\bigr]\sinh(\vartheta\phi),
 \qquad j_K=2KI_\omega(r)\cosh(\vartheta\phi)\sinh(qA).
 \label{eq:rotor-full-moments}
\end{equation}
Testing the electric equation by \(\phi\) forces \(\phi=0\) for
every equilibrium.  Only then does the magnetic equation reduce to
\begin{equation}
 L_AA=2KI_\omega(r)\sinh(qA).                  \label{eq:rotor-scalar}
\end{equation}
Put
\begin{equation}
 w_\omega=2qI_\omega,\quad
 c_\omega=\tfrac13q^3I_\omega=\tfrac16q^2w_\omega,\quad
 K_*=\Lambda_1(w_\omega).                         \label{eq:rotor-threshold}
\end{equation}
The positive eigenfunction \(L_A\psi=K_*w_\omega\psi\)
gives crossing form \(-\int w_\omega\psi^2<0\) and
\begin{equation}
 A_\pm(s)=\pm s\psi+O(s^3),\qquad
 K(s)=K_*-K_*
 \frac{\int_\Omega c_\omega\psi^4\dd x}
      {\int_\Omega w_\omega\psi^2\dd x}s^2+O(s^4).
 \label{eq:rotor-branch-direction}
\end{equation}
Thus the pitchfork is subcritical.  The profile satisfies
Proposition \ref{prop:superlinear-case-M}, since
\(F(r,s)/(w_\omega s)=\sinh(qs)/(qs)>1\).  Its positive
and negative global continua have unbounded magnetic-potential
norm but remain in \(0<K<K_*\); \(\phi=E=0\) identically.
In particular, with \(q_+=\omega r_+\),
\begin{equation}
 M_R=\frac{\sinh(q_+R)}{q_+R},\qquad
 0<K<K_*\frac{q_+R}{\sinh(q_+R)}
 \ \Longrightarrow\ \|A\|_\infty>R\quad(A>0).
 \label{eq:rotor-radius-test}
\end{equation}
The parameter need not approach zero along the unbounded end.
On a circular torus, Proposition \ref{prop:general-torus}
provides Bessel-zero enclosures for \(K_*\).
No growing RVM mode follows from this stationary calculation:
the orbit-projected dynamical operator on the solid torus is not
evaluated here.
\end{example}

\subsection{Case EM: simultaneous electric and magnetic norm escape}

The preceding magnetic class can be coupled to a growing charge
imbalance without losing positivity or computability.  Fix a desired
bifurcation parameter \(K_*>0\) and constants
\begin{equation}
 \sigma>\gamma>0,
 \qquad
 \delta(K)=\gamma(K-K_*),\qquad
 \beta(K)=\Lambda e^{\sigma(K-K_*)},\qquad
 a_\pm(K)=e^{\pm\delta(K)}.                       \label{eq:double-escape-parameters}
\end{equation}
Use $C$ and the integrable, strictly decreasing curvature $g=G''$
from the sublinear coefficient-growth class
\eqref{eq:unbounded-G}--\eqref{eq:unbounded-profiles}, not the
positive-tail or superlinear thermal classes.  Define
\begin{equation}
 \mu_K^\pm(e,p)=a_\pm(K)e^{-e}
                   \bigl(C+\beta(K)G(p)\bigr).    \label{eq:double-escape-profiles}
\end{equation}
These are strictly positive, energy-decreasing two-species profiles.
Set
\begin{equation}
 N_K(r,s)=\int_{\R^3}e^{-\av v}
       \bigl(C+\beta(K)G(r(v_\varphi+s))\bigr)\dd v.
 \label{eq:double-escape-N}
\end{equation}
The function \(N_K(r,\cdot)\) is even, and
\(N_K(r,s)\ge Cm_0\).  The identities used to derive
\eqref{eq:unbounded-scalar} now give the exact coupled system
\begin{equation}
 \left\{
 \begin{aligned}
 L_\phi\phi
  &=2N_K(r,A)\sinh\bigl(\delta(K)-\phi\bigr),\\
 L_AA
  &=\beta(K)\cosh\bigl(\delta(K)-\phi\bigr)H(r,A),\\
 A|_{\partial\Omega}&=\phi|_{\partial\Omega}=0.
 \end{aligned}\right.                            \label{eq:double-escape-system}
\end{equation}
For each \((A,K)\), the first equation has a unique solution
\(\phi=\Phi_K[A]\).  Indeed its field derivative is
\[
 L_\phi+2N_K(r,A)\cosh(\delta(K)-\phi)>0.
\]
Write \(\phi_0(K)=\Phi_K[0]\); then
\((0,\phi_0(K))\) is the reference branch and
\(\phi_0(K_*)=0\).

\begin{theorem}[A coupled curve with both potentials unbounded]
\label{thm:coupled-double-escape}
For the profiles \eqref{eq:double-escape-profiles}, the reference
equilibrium \((0,0,K_*)\) is a simple transverse magnetic bifurcation
point.  In the cone \(A>0\), the right-hand local branch belongs to a
global connected component \(\mathscr C_+^{EM}\) whose parameter
projection is unbounded above.  Every point of this component with
\(K\ge K_*\) satisfies
\begin{equation}
 0\le\phi\le\delta(K).                            \label{eq:double-escape-phi-barrier}
\end{equation}
For every \(R>0\), the following completely computable implications
hold on \(\mathscr C_+^{EM}\):
\begin{align}
 2Cm_0\sinh\bigl(\delta(K)-R\bigr)
   &>\lambda_{\phi,1}R
 &&\Longrightarrow\quad \|\phi\|_\infty>R,
 \label{eq:double-electric-test}\\
 \beta(K)&>\frac{\Lambda}{q_R}
 &&\Longrightarrow\quad \|A\|_\infty>R.
 \label{eq:double-magnetic-test}
\end{align}
Consequently, along every sequence in \(\mathscr C_+^{EM}\) with
\(K\to\infty\),
\begin{equation}
 \|\phi\|_{L^\infty}+\|\phi\|_{X_\phi}\to\infty,
 \qquad
 \|A\|_{L^\infty}+\|A\|_{X_A}\to\infty.        \label{eq:double-escape-conclusion}
\end{equation}
Thus this component exhibits simultaneous electric and magnetic norm
escape.
\end{theorem}

\begin{proof}
At \(K=K_*\), the activities agree, \(\phi_0=0\), and the mixed
linearized blocks vanish by parity.  The electric block is strictly
positive, while the magnetic block is
\(L_A-\Lambda\alpha\).  Its kernel is generated by the positive
weighted principal eigenfunction \(\psi\).  Along the reference branch,
the derivative of the magnetic coefficient is
\[
 \left.\frac{\dd}{\dd K}\bigl[
 \beta(K)\cosh(\delta(K)-\phi_0(K))\bigr]
 \right|_{K=K_*}=\sigma\Lambda;
\]
the derivative of the cosh factor vanishes at zero.  Hence the crossing
form is
\(-\sigma\Lambda\int_\Omega\alpha\psi^2\dd x<0\).

Solving the electric equation as \(\phi=\Phi_K[A]\) reduces the system
to a compact strongly positive magnetic equation on $Y_+$; here
$\beta(K)>0$ throughout the family.  The
Dirichlet resolvent maps bounded subsets compactly into $Y$;
the strong maximum principle and Hopf lemma put every nonzero
nonnegative solution in the cone interior.  For \(K>K_*\), its
linearization at zero has negative principal eigenvalue because
\(\beta(K)>\Lambda\) and
\(\cosh(\delta-\phi_0)\ge1\).  On large balls the map is compressive:
\(H(r,s)\) is uniformly bounded in \(s\), while
\(0\le\Phi_K[A]\le\delta(K)\) and hence the cosh factor is at most
\(\cosh\delta(K)\).  The positive-cone theorem
\cite{DancerPositive1973,Amann} therefore gives a right-hand positive
continuum issuing from the local germ.

There are no positive solutions for \(K<K_*\).  Indeed, the maximum
principle gives \(\delta(K)\le\phi\le0\), so
\[
 \cosh(\delta(K)-\phi)\le e^{\gamma(K_*-K)}.
\]
Since \(H(r,A)\le\alpha(r)A\) for \(A>0\), pairing the magnetic
equation with \(\psi\) would give
\[
 \Lambda\int_\Omega\alpha\psi A\dd x
 \le \Lambda e^{-(\sigma-\gamma)(K_*-K)}
       \int_\Omega\alpha\psi A\dd x,
\]
which is impossible because \(\sigma>\gamma\).  At \(K=K_*\), the
strict inequality \(H(r,A)<\alpha(r)A\) rules out a nonzero positive
solution.  At any $K_1>K_*$ the reference magnetic derivative is
\[
 L_A-\beta(K_1)\cosh\bigl(\delta(K_1)-\phi_0(K_1)\bigr)\alpha(r).
\]
Its principal eigenvalue is strictly negative: $\beta(K_1)>\Lambda$
and the cosh factor is at least one.  If positive solutions approached
$(0,\phi_0(K_1),K_1)$, division of the magnetic equations by their
$C^1$ norms and elliptic compactness would give a nonzero nonnegative
kernel vector for this derivative.  Positivity improvement would make
it a principal eigenfunction, contradicting the negative principal
eigenvalue.  This also excludes return at a higher-mode zero
eigenvalue.  Thus the positive component cannot meet another reference
state or the cone boundary.  No admissibility boundary is present.  On every bounded
\(K\)-interval, the bound on \(H\), the maximum-principle bound on
\(\phi\), the at-most-linear growth of \(G\), and elliptic estimates
bound both fields in \(\calX\).  The unilateral global alternative
therefore forces the parameter projection of
\(\mathscr C_+^{EM}\) to be unbounded above.

It remains to identify which components escape.  To prove
\eqref{eq:double-electric-test}, assume its strict hypothesis and
\(\|\phi\|_\infty\le R\).  The strict hypothesis implies
\(\delta(K)>R\); hence the electric equation and \(N_K\ge Cm_0\),
paired with the positive first Dirichlet eigenfunction \(\varphi_1\),
give
\[
 \lambda_{\phi,1}R\int_\Omega\varphi_1\dd x
 \ge \int_\Omega\varphi_1\rho_K\dd x
 \ge2Cm_0\sinh(\delta(K)-R)
       \int_\Omega\varphi_1\dd x,
\]
a contradiction.  Similarly, if \(0<A\le R\), then
\(H(r,A)\ge q_R\alpha(r)A\) and \(\cosh(\delta-\phi)\ge1\).
Pairing the magnetic equation with the positive eigenfunction
\(L_A\psi=\Lambda\alpha\psi\) yields
\[
 \Lambda\int_\Omega\alpha\psi A\dd x
 =\beta(K)\int_\Omega\psi\cosh(\delta-\phi)H(r,A)\dd x
 \ge\beta(K)q_R\int_\Omega\alpha\psi A\dd x,
\]
so \(\beta(K)\le\Lambda/q_R\), proving
\eqref{eq:double-magnetic-test}.  Both thresholds are eventually
satisfied because \(\delta(K)\to\infty\) and
\(\beta(K)\to\infty\).  The embeddings into \(C^0\) complete
\eqref{eq:double-escape-conclusion}.
\end{proof}

The family \eqref{eq:double-escape-profiles} combines the electric-only
J\"uttner mechanism with the sublinear magnetic class: setting \(G=0\)
leaves the unequal-activity electric graph, while setting
\(\gamma=0\) restores the neutral magnetic class.  With both effects
present, the activity ratio
\(a_+/a_-=e^{2\delta}\) drives charge separation and the coefficient
\(\beta\) drives canonical-momentum current.  The inequalities
\eqref{eq:double-electric-test} and
\eqref{eq:double-magnetic-test} show separately, rather than through a
combined norm, that neither field can remain bounded.

\end{document}